\documentclass[11pt]{article}

\newcommand{\bb}{\textbf}

\newcommand{\ii}{\textit}

\newcommand{\hh}{\hspace}
\newcommand{\vv}{\vspace}

\usepackage[paper=a4paper, margin=2cm]{geometry}

\usepackage{graphicx}
\usepackage[colorlinks=true, linkcolor=blue, urlcolor=blue, citecolor = blue]{hyperref}
\usepackage{titlesec}
\usepackage{textcomp}
\usepackage{color}
\usepackage{colortbl}
\usepackage[dvipsnames]{xcolor}
\usepackage {amssymb, amsthm,amsmath}
\usepackage{bbding}
\usepackage{enumerate}
\usepackage{esint}

\usepackage{tikz}
\usepackage{pifont}
\usepackage{wasysym}

\makeatletter
\newcommand{\thickhline}{%
    \noalign {\ifnum 0=`}\fi \hrule height 1pt
    \futurelet \reserved@a \@xhline
}
\newcolumntype{"}{@{\hskip\tabcolsep\vrule width 1pt\hskip\tabcolsep}}
\makeatother

\begin{document}
\setlength{\parskip}{0pt}
% \raggedbottom

\title{Mathematical model for delayed responses in immune checkpoint blockades}

\author{Collin Y. Zheng$^{1}$         \and
        Peter S. Kim$^{2}$
}

\maketitle

\noindent $^{1}$ \quad School of Mathematics and Statistics, The University of Sydney, NSW 2006, Australia, \\ (collin.zheng@sydney.edu.au), \\
\noindent $^{2}$ \quad School of Mathematics and Statistics, The University of Sydney, NSW 2006, Australia, \\ (peter.kim@sydney.edu.au).

\begin{abstract}
We introduce a set of ordinary differential equations (ODE) that qualitatively reproduces delayed responses observed in immune checkpoint blockade therapy (e.g. anti-CTLA-4 Ipilimumab). This type of immunotherapy has been at the forefront of novel and promising cancer treatments over the past decade and was recognised by the 2018 Nobel Prize in Medicine. Our model describes the competition between effector T cells and non-effector T cells in a tumour. By calibrating a small subset of parameters that control immune checkpoint expression along with the patient’s immune-system cancer readiness, our model is able to simulate either a complete absence of patient response to treatment, a quick anti-tumour T cell response (within days) or a delayed response (within months). Notably, the parameter space that generates a delayed response is thin and must be carefully calibrated, reflecting the observation that a small subset of patients experience such reactions to checkpoint blockade therapies. Finally, simulations predict that the anti-tumour T cell storm that breaks the delay is very short-lived compared to the length of time the cancer is able to stay suppressed. This suggests the tumour may subsist off an environment hostile to effector T cells; however, these cells are---at rare times---able to break through the tumour immunosuppressive defences to neutralise the tumour for a prolonged period. Our simulations aim to qualitatively describe the delayed response phenomenon without making precise fits to particular datasets, which are limited. It is our hope that our foundational model will stimulate further interest within the immunology modelling field.
\end{abstract}


\noindent\textbf{Acknowledgements:} This scientific work was supported by the Australian Postgraduate Award (CYZ) and the Australian Research Council Discovery Project DP180101512 (PSK). We thank Peter P. Lee (Beckman Research Institute, City of Hope, California, USA) for his guidance and insights on immune checkpoint blockades.

%%%%%%%

\section{Introduction}

% Cancer
Cancer has been a leading cause of death around the world for decades, yet much remains unknown about its mechanisms of establishment and destruction. While traditional modes of treatment---surgery, chemotherapy and radiotherapy---have played important roles in treatment, they often do not result in a definitive cure. Indeed, in a substantial number of cases, patients who experience tumour regression are later met with news that cancer has returned at the original site or elsewhere through metastasis. The prognosis for patients in the latter category remains low for all cancer types~\cite{dillekasAre90Deaths2019}.\\

% Cancer and the immune system
An important piece of the puzzle surrounding cancer concerns our immune system: why does our anti-tumour immunity appear so incompetent? A clue lies in our immune system's fundamental design principle of immune surveillance, the notion that the entire system is programmed to only detect and destroy non-self cells. Yet it is precisely the notion that cancer does not belong in that category that causes such a dilemma. The inconvenient truth is cancer cells are not foreign, but are self cells that have gone rogue. Seen through this prism, the frequent failure of anti-tumour immunity to prevent prevent tumour escape should not come as a surprise.~\cite{swannImmuneSurveillanceTumors2007}\\

% Focus on T cells
With that said, the immune system has shown great potential. Three types of effector white blood cells are known to be capable of destroying cancer cells: natural killer cells (NK cells), macrophages (Type 2 phenotype macrophages) and T cells (CD8$^+$ T-lymphocytes). Out of these three, it is the T cell that shows the most promise~\cite{gonzalezRolesImmuneSystem2018, waldmanGuideCancerImmunotherapy2020}.\\
 
% Cancer hiding from T cells
Cancer employs a wide array of strategies to evade and hide from T cells, ultimately resulting in tumour escape. These include its weak antigenicity that permits it to mask its identity, tendency for rapid mutation that makes it a difficult moving target, and its ability to leverage a wide array of immunosuppressive mechanisms to deal with the T cells that have found the tumour site~\cite{nurievaTcellToleranceCancer2013}.\\

% Immunotherapy
Over the past decade, a new brand of treatment---immunotherapy---has aimed to reverse the immune system's bad fortunes by placing it at the fore in the struggle against cancer~\cite{topalianCancerImmunotherapyComes2011, alexanderCheckpointImmunotherapyRevolution2016b, waldmanGuideCancerImmunotherapy2020}. Immunotherapy strives to understand the dynamic interplay between the immune system and cancer so that novel treatments centered around stimulating robust anti-tumour immune responses can be developed. Immunotherapeutic treatments include vaccines, cytokines and antibodies that modulate an adaptive anti-tumour immune response. The latter, known as immune checkpoint inhibitor therapy, which aims to shift the balance of power back towards T cells, is the focus of current research efforts and has shown tremendous potential. Research and trialling in this space has exploded since 2010, culminating in the 2018 Nobel Prize in Medicine.\\

% Intro to checkpoint inhibitors
T cells express special receptor molecules on their surface called immune checkpoints, which serve to suppress their function and proliferation. The two most notable immune checkpoints are cytotoxic T-lymphocyte-associated protein 4 (CTLA-4) and programmed cell death protein 1 (PD-1). The purpose of these two receptors is to assist in shutting down the immune system once an infection has passed. This is an important feature of the immune system as white blood cells that live past their due date become toxic. Unfortunately, cancer cells have the ability to bind with checkpoint receptors like PD-1, giving them an easy way to deny an anti-tumour T-cell response~\cite{pardollBlockadeImmuneCheckpoints2012, sliwkowskiAntibodyTherapeuticsCancer2013, marzoliniCHECKPOINTBLOCKADECANCER2015, callahanCTLA4PD1Pathway2015, buchbinderCTLA4PD1Pathways2016a, seidelAntiPD1AntiCTLA4Therapies2018a}. In essence, the cancer cripples its adversary before the battle has begun.\\

\begin{figure}
    \centering
    \includegraphics[scale=1]{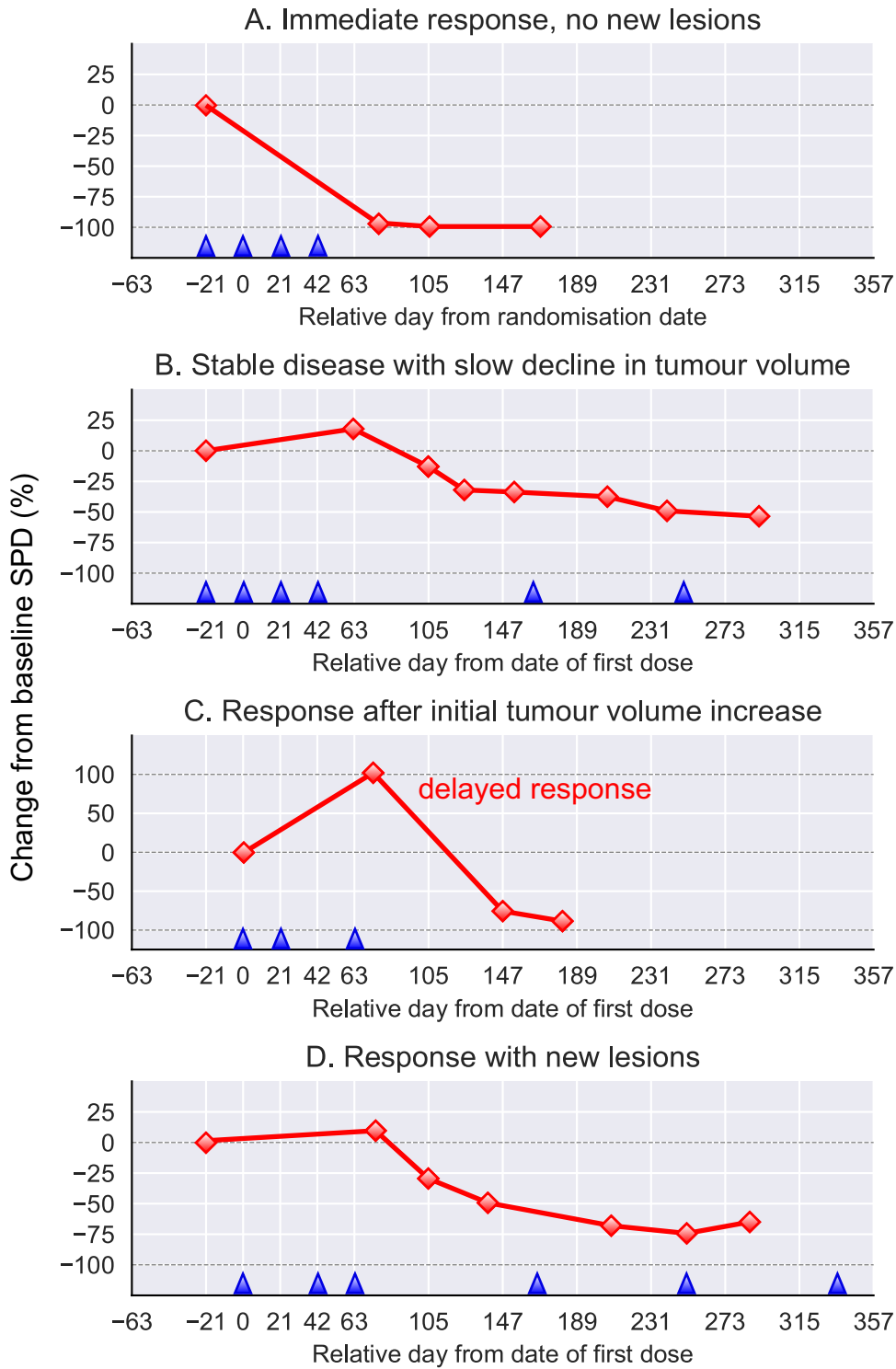}
    \caption{\bb{Responses to the CTLA-4 inhibitor Ipilimumab.} Trials from \cite{wolchokGuidelinesEvaluationImmune2009} showed four distinct patterns of clinical response from anti-CTLA-4 Ipilimumab trial data: A) immediate response with no new lesions; B) stable disease with slow decline in tumor volume; C) response after initial tumour increase; and D) response with new lesions. The modelling in this paper is motivated by the delayed response seen in case C. The blue triangles indicate the administration times of Ipilimumab (3-4 dosages every 3 weeks, sometimes followed by maintenance doses). Note: the SPD (Sum of Product of Diameters), used in the \ii{response assessment in neuro-oncology (RANO) criteria} for assessing cancer treatment efficacy, is calculated as the sum of the products of maximal diameters of lesions in the patient.}
    \label{SpikeODE_intro_Wolchok}
\end{figure}

% Checkpoint inhibitors in detail
Checkpoint inhibitor therapy targets T cells in both lymph nodes and the tumour site by blocking CTLA-4 or PD-1 with monoclonal antibodies. Ipilimumab, which blocks CTLA-4, was approved by the United States Food and Drug Administration (FDA) in 2011, while the anti-PD-1 antibodies Nivolumab and Pembrolizumab were FDA-approved in 2014~\cite{alexanderCheckpointImmunotherapyRevolution2016a}. Alongside these and several other approved antibodies, a wide range of new ones are currently being developed and trialled around the world. Although durable outcomes have been observed in many checkpoint inhibitor trials, a delayed response has been observed in a small number of cases, where a tumour continues to grow for a period of time before the antibody treatment appears to take effect\cite{wolchokGuidelinesEvaluationImmune2009} (see Fig. \ref{SpikeODE_intro_Wolchok}). The mechanisms behind this mysterious lag remain a mystery and efforts of modelling in the area are nascent.\\

% Scope of paper
In this paper, we introduce a system of ordinary differential equations (ODE) that describes tumour-T-cell dynamics and incorporates CTLA-4 and PD-1 inhibitors. Specifically, we are interested in motivating a relatively simple set of equations that is capable of qualitatively reproducing delayed responses from checkpoint inhibitors. Rather than make a precise fit to specific data sets, which is limited in experiments, we aim to produce a set of characteristic---or \ii{idealised}---simulations that capture the qualitative nature of experimental data without precisely fitting to them. It is our hope that our foundational model can pique the interest of mathematical biologists in a nascent space in the cancer modelling field. Our model predicts a thin region in the parameter space that admits the delayed-response case. Since parameter combinations can be thought of as mapping to specific immune profiles, our results suggest it may be relatively rare for patients to exhibit delayed responses and moreover, it may be difficult or not possible for an individual that starts sufficiently far away from this threshold parameter region to respond to immune checkpoint blockade therapy. Moreover, the model simulates T cell dynamics in the tumour reminiscent of the competition between effector T cells and regulatory T cells~\cite{oleinikaSuppressionSubversionEscape2013, quezadaCTLA4BlockadeGMCSF2006, curranTumorVaccinesExpressing2009}.

% %%%%%%%%%%%%%%%%%%%%%%%%%%%%%%%%%%%%%%%%%%%%%%%%%%%%%%%%%%%%%%%%%%%%%%%%%%%%%
% Intro DIAGRAM: Wolchonk delays
% %%%%%%%%%%%%%%%%%%%%%%%%%%%%%%%%%%%%%%%%%%%%%%%%%%%%%%%%%%%%%%%%%%%%%%%%%%%%%

%%%%%%%%%%%%%%%%%%%%%%%%%%%%%%%%%%%%%%%%%%%%%%%%%%%%%%%%%%%%%%%%%
%%%%%%%%%%%%%%%%%%%%%%%%%%%%%%%%%%%%%%%%%%%%%%%%%%%%%%%%%%%%%%%%%
%%%%%%%%%%%%%%%%%%%%%%%%%%%%%%%%%%%%%%%%%%%%%%%%%%%%%%%%%%%%%%%%%
%%%%%%%%%%%%%%%%%%%%%%%%%%%%%%%%%%%%%%%%%%%%%%%%%%%%%%%%%%%%%%%%%
%%%%%%%%%%%%%%%%%%%%%%%%%%%%%%%%%%%%%%%%%%%%%%%%%%%%%%%%%%%%%%%%%
%%%%%%%%% 2. MODEL %%%%%%%%%%%%%%%%%%%%%%%%%%%%%%%%%%%%%%%%%%%%%%
%%%%%%%%%%%%%%%%%%%%%%%%%%%%%%%%%%%%%%%%%%%%%%%%%%%%%%%%%%%%%%%%%
%%%%%%%%%%%%%%%%%%%%%%%%%%%%%%%%%%%%%%%%%%%%%%%%%%%%%%%%%%%%%%%%%
%%%%%%%%%%%%%%%%%%%%%%%%%%%%%%%%%%%%%%%%%%%%%%%%%%%%%%%%%%%%%%%%%
%%%%%%%%%%%%%%%%%%%%%%%%%%%%%%%%%%%%%%%%%%%%%%%%%%%%%%%%%%%%%%%%%
%%%%%%%%%%%%%%%%%%%%%%%%%%%%%%%%%%%%%%%%%%%%%%%%%%%%%%%%%%%%%%%%%

\newpage
\section{Model}\label{sec:model}

\begin{figure}
    \centering
    \includegraphics[scale=1]{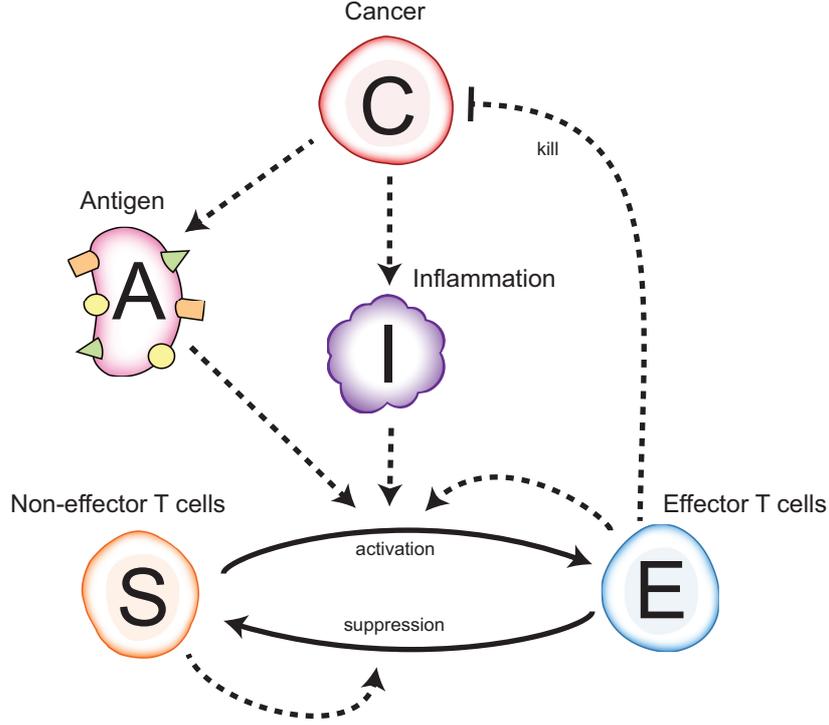}
    \caption{\bb{Diagram of the model}. The dotted line represent a cause or effect, while the solid lines represent the transition of one T cell type to another. In our model, cancer ($C$) causes the upregulation of tumour-specific antigen ($A$) and inflammation ($I$), which affects the population dynamics between effector ($E$) and non-effector ($S$) T cells. Effector T cells kill $C$ cells, causing downward pressure on its population.}
    \label{SpikeODE_model_model}
\end{figure}

We propose a mathematical model to help understand the dynamics of the delayed immune responses observed in some immune checkpoint inhibitor therapy trials~\cite{wolchokGuidelinesEvaluationImmune2009}. We present the model as a set of ODEs with five variables. These variables are summarised in Table~\ref{table:variables}. In short, the model describes \bb{cancer} ($C$) giving rise to antigen and inflammatory \bb{signals} ($A$ and $I$) that affect the interaction and competition between effector and non-effector \bb{T cells} ($E$ and $S$). A diagram of the model can be found in Fig. \ref{SpikeODE_model_model}.

% %%%%%%%%%%%%%%%%%%%%%%%%%%%%%%%%%%%%%%%%%%%%%%%%%%%%%%%%%%%%%%%%%%%%%%%%%%%%
% Model DIAGRAM: Model diagram
% %%%%%%%%%%%%%%%%%%%%%%%%%%%%%%%%%%%%%%%%%%%%%%%%%%%%%%%%%%%%%%%%%%%%%%%%%%%%

\vv{2pc}

% %%%%%%%%%%%%%%%%%%%%%%%%%%%%%%%%%%%%%%%%%%%%%%%%%%%%%%%%%%%%%%%%%%%%%%%%%%%%
% Model EQUATIONS: The model
% %%%%%%%%%%%%%%%%%%%%%%%%%%%%%%%%%%%%%%%%%%%%%%%%%%%%%%%%%%%%%%%%%%%%%%%%%%%%

\begin{align}
\frac{dC}{dt} &= f(C)\,C - \kappa\,CE, \label{eq:spikeC}\\
\frac{dA}{dt} &= r_A\,C - \delta_A\,A, \label{eq:spikeA}\\
\frac{dI}{dt} &= r_I\,CE - \delta_I\,I, \label{eq:spikeI}\\
\frac{dE}{dt} &= -r_E\,(E - E^*) + \beta AIES - \gamma ES, \label{eq:spikeE}\\
\frac{dS}{dt} &= -r_S\,(S - S^*) - \beta AIES + \gamma ES. \label{eq:spikeS}
\end{align}

\vv{1.5pc}

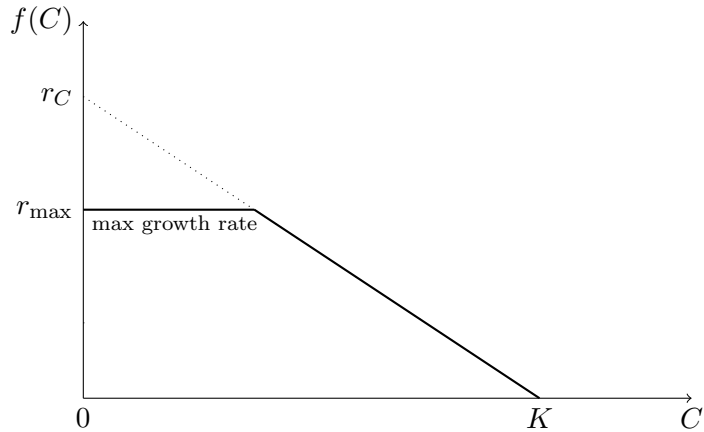
\begin{figure}[htbp]
    \centering
    \begin{tikzpicture}
    
    % axes
    \draw[->] (0,0) -- (8,0) node[anchor=north] {$C$};
    \draw[->] (0,0) -- (0,5) node[anchor=east] {$f(C)$};

    % labels
    \draw	(0,0) node[anchor=north] {0};
    \draw	(6,0) node[anchor=north] {$K$};
    \draw	(0,4) node[anchor=east] {$r_C$};
    \draw	(0,2.5) node[anchor=east] {$r_{\text{max}}$};

    % words
    \draw   (1.2,2.3) node{{\scriptsize max growth rate}};
    
    % line
    \draw[dotted] (0,4) -- (6,0);
    \draw[dotted] (0,1) -- (0,1);
    \draw[thick] (0,2.5) -- (2.25,2.5);
    \draw[thick] (2.25,2.5) -- (6,0);
    
    % % Us
    % \draw[thick] (0,0) -- (2,2) -- (6,2);
    % \draw (1,1.5) node {$U_s$}; %label
    
    % % Psis
    % \draw[thick,dashed] (0,3) -- (2,3) parabola[bend at end] (6,1);
    % \draw (2.5,3) node {$\varPsi_s$}; %label
    
    \end{tikzpicture}
    \caption{Plot of the cancer growth rate $\displaystyle f(C) = \text{min}\left\{r_C\left(1 - \frac{C}{C^*}\right),\,\, r_\text{max}\right\}$.}
    \label{SpikeODE_pic_RGR}
\end{figure}

%%% EQN 1 CANCER
\indent \bb{Equation \ref{eq:spikeC}} describes the concentration of the cancer population, $C$, in the tumour microenvironment. For simplicity, we assume this growth occurs in the form of a solid tumour mass. The first term represents the intrinsic growth of the tumour. We assume this growth to be a function that saturates at $r_{max}$, and set the growth rate to be $$f(C) = \text{min}\left\{r_C\left(1 - \frac{C}{C^*}\right),\,\,r_\text{max}\right\}.$$ We elaborate on this choice for $f(C)$ in Section 5. A plot of $f(C)$ versus $C$ is shown in Figure~\ref{SpikeODE_pic_RGR}. The second term represents the mass-action rate at which the cancer is killed by effector T cells with kinetic coefficient $\kappa$.

%%% EQN 2 ANTIGEN PRESENTATION (pMHC)
\indent \bb{Equation \ref{eq:spikeA}} describes the concentration of cancer-specific antigen, $A$, that is presented by antigen-presenting cells (APC). The antigen presented by APCs---particularly dendritic cells---serve to activate naive anti-cancer T cells, which hunt and kill cancer cells in the tumour microenvironment~\cite{paluckaCancerImmunotherapyDendritic2012}. In the first term, we assume cancer-associated antigen enters the system at a rate proportional to the cancer population with coefficient $r_A$. The second term describes antigen degrading at rate $\delta_A$.

%%% EQN 3 INFLAMMATION
\indent \bb{Equation \ref{eq:spikeI}} describes the volume of inflammatory signals, $I$, at the tumour site. These signals activate immature dendritic cells (DC), resulting in their migration to the lymph nodes in search of naive T cells to activate~\cite{tangHighmobilityGroupBox2010, simsHMGB1RAGEInflammation2010}. The first term represents the rate inflammatory signal is generated due to the killing of cancer cells by effector T cells (Teff).  Like cancer killing, we assume this signal is produced at a mass-action rate with coefficient $r_I$. The second term describes the inflammation degrading at rate $\delta_I$.

%%% EQN 4 and 5 T CELLS %%
\indent \bb{Equations \ref{eq:spikeE} and \ref{eq:spikeS}} model the effector (immunogenic) and non-effector (tolerogenic) T cell concentrations, $E$ and $S$, respectively. For simplicity we regard $E$ and $S$ as two different \ii{states} for T cells, where effector T cells have the capacity to kill cancer cells, non-effector T cells do not, and T cells can transition from one state to the other. In general, studies have suggested tumour microenvironments play host to complex interactions between immunogenic and tolerogenic immune cells.\cite{balkwillTumorMicroenvironmentGlance2012, hornyakRoleIndoleamine23Dioxygenase2018, tesiMDSCMostImportant2019}. For instance, CD8$^+$ T cells and regulatory T cells are important respective examples of immunogenic and tolerogenic T cells~\cite{grauerCD4FoxP3Regulatory2007, oleinikaSuppressionSubversionEscape2013}. In the first terms for Equations \ref{eq:spikeE} and \ref{eq:spikeS}, we assume both T cell populations hold a tendency towards some sort of immune homeostasis. In the case of $E$, we assume the immunogenic cancer mobilises a base level T cell response, so that the population approaches some steady state level $E^*$ at a rate $r_E$. Similarly, we assume the immune system generates a base level of tolerogenic T cells that approach some steady state $S^*$ at a rate $r_S$. The second terms in Equations \ref{eq:spikeE} and \ref{eq:spikeS} describe T cells switching states from tolerogenic to immunogenic, i.e. transitioning from $S$ to $E$. Here, we assume that $S$ cells are activated by a cocktail of signals, $g(A, I, E)$, consisting of antigen, $A$; inflammation, $I$; and positive growth signals, assumed to be proportional to $E$. For simplicity, we take this rate to be proportional to the product of signals, so that the full term is of the form $g(A,I,E)\,S = \beta\,AIES$ for some constant $\beta$.  Finally, the third terms in Equations \ref{eq:spikeE} and \ref{eq:spikeS} describe T cells switching states from immunogenic to tolerogenic, i.e., T cells transitioning from $E$ to $S$. For simplicity, we assume that $E$ cells are being anergised or suppressed by inhibitory signals proportional to $S$, so that the full term takes the form $\gamma\,SE$ for some constant $\gamma$.

%%% CTLA-4 and PD-1
\indent The roles played by the immune checkpoint receptors CTLA-4 and PD-1 are described by the coefficients $\beta$ and $\gamma$, respectively. It is believed that CTLA-4 primarily plays a role in the activation of effector T cells, while PD-1 plays part in their anergisation and suppression~\cite{pardollBlockadeImmuneCheckpoints2012, postowImmuneCheckpointBlockade2015, buchbinderCTLA4PD1Pathways2016a, seidelAntiPD1AntiCTLA4Therapies2018a}. Hence the effect of CTLA-4 and their antibody inhibitors (e.g., Ipilimumab) is calibrated by the value of $\beta$, which controls the activation rate of effector T cells, while the effect of PD-1 and their inhibitors (e.g., Nivolumab) is calibrated by the value of $\gamma$, which controls the suppression rate of effector T cells.

% %%%%%%%%%%%%%%%%%%%%%%%%%%%%%%%%%%%%%%%%%%%%%%%%%%%%%%%%%%%%%%%%%%%%%%%%%%%%
% Model DIAGRAM: f(C)
% %%%%%%%%%%%%%%%%%%%%%%%%%%%%%%%%%%%%%%%%%%%%%%%%%%%%%%%%%%%%%%%%%%%%%%%%%%%%

% %%%%%%%%%%%%%%%%%%%%%%%%%%%%%%%%%%%%%%%%%%%%%%%%%%%%%%%%%%%%%%%%%%%%%%%%%%%%
% Model VARIABLES
% %%%%%%%%%%%%%%%%%%%%%%%%%%%%%%%%%%%%%%%%%%%%%%%%%%%%%%%%%%%%%%%%%%%%%%%%%%%%

\begin{table}[h]
\centering
    \begin{tabular}{l l l}
    \hline\noalign{\smallskip} 
    Variable & Description & Units\\
    \noalign{\smallskip}\thickhline\noalign{\smallskip}
    $C(t)$ & Cancer & cells/nL \smallskip\\
    $A(t)$ & Antigen presentation & peptides/nL \smallskip\\
    $I(t)$ & Inflammation & ng/nL \smallskip\\
    $E(t)$ & Effector (immunogenic) T cells & cells/nL \smallskip\\
    $S(t)$ & Non-effector (tolerogenic) T cells & cells/nL \smallskip\\
    \hline
  \end{tabular}\smallskip
  \caption{\bb{Summary of variables}.}
  \label{table:variables}
\end{table}

%%%%%%%%%%%%%%%%%%%%%%%%%%%%%%%%%%%%%%%%%%%%%%%%%%%%%%%%%%%%%%
%%%%%%%%%%%%%%%%%%%%%%%%%%%%%%%%%%%%%%%%%%%%%%%%%%%%%%%%%%%%%%
%%%%%%%%%%%%%%%%%%%%%%%%%%%%%%%%%%%%%%%%%%%%%%%%%%%%%%%%%%%%%%
%%%%%%%%%%%%%%%%%%%%%%%%%%%%%%%%%%%%%%%%%%%%%%%%%%%%%%%%%%%%%%
%%%%%%%%%%%%%%%%%%%%%%%%%%%%%%%%%%%%%%%%%%%%%%%%%%%%%%%%%%%%%%
%%%%%%%%% 3. PARAMETER ESTIMATES %%%%%%%%%%%%%%%%%%%%%%%%%%%%%
%%%%%%%%%%%%%%%%%%%%%%%%%%%%%%%%%%%%%%%%%%%%%%%%%%%%%%%%%%%%%%
%%%%%%%%%%%%%%%%%%%%%%%%%%%%%%%%%%%%%%%%%%%%%%%%%%%%%%%%%%%%%%
%%%%%%%%%%%%%%%%%%%%%%%%%%%%%%%%%%%%%%%%%%%%%%%%%%%%%%%%%%%%%%
%%%%%%%%%%%%%%%%%%%%%%%%%%%%%%%%%%%%%%%%%%%%%%%%%%%%%%%%%%%%%%%
%%%%%%%%%%%%%%%%%%%%%%%%%%%%%%%%%%%%%%%%%%%%%%%%%%%%%%%%%%%%%%%

\section{Parameter estimates}\label{sec:parameters}

The majority of parameters for $C$, $A$ and $I$ were calculated from biological data, with details below. The remaining parameters, notably $\beta$ (associated with T cell activation and CTLA-4) and $\gamma$ (associated with effector T cell suppression and PD-1) were fitted using MATLAB's non-linear least squares tool, \ii{lsqnonlin}. Since we are particularly interested in the delayed-response case, we fit our simulations around a two-month delay representing a combination Ipilimumab (anti-CTLA-4) and Nivolumab (anti-PD-1) treatment. For the rest of this paper, we designate the final fitted set of parameters giving rise to this two-month delay the \bb{baseline parameters}. These are summarised in Table \ref{table:baselineparams}. Moreover, a summary of the initial values can be found in Table \ref{table:ICs}. As we shall see in the Results and Discussion sections, all simulations for all response types may be obtained by perturbing this baseline set of parameters.

%%%%%%%%%%%%%%%%%%%%%%%%%%%%%%%%%%%%%%%%%%%%%%%%%%%%%%%%%%%%%%
% Param TABLE: baseline parameter set that generates 2 month delay
%%%%%%%%%%%%%%%%%%%%%%%%%%%%%%%%%%%%%%%%%%%%%%%%%%%%%%%%%%%%%%

\begin{table}[htbp]
\centering
  \begin{tabular}{l p{6cm} l l l}
  
    \hline\noalign{\smallskip} 
    % TUMOUR
    Symbol & Description & Estimate\\
    \noalign{\smallskip}\thickhline\noalign{\smallskip}
    $r_C$ & logistic growth rate of cancer & 1.0 day$^{-1}$ \smallskip\\
    $r_\text{max}$ & maximum growth rate of cancer & 0.09 day$^{-1}$\smallskip\\ % owenMathematicalModellingUse2004
    $C^*$ & cancer steady state concentration & $10^3$ cells/nL\smallskip\\
    $\kappa$ & killing rate of cancer cells by T cells & 1.2 nL/(cells $\cdot$ day)\smallskip\\
    
    \hline\noalign{\smallskip}
    % antigen presentation
    $r_A$ & antigen presentation source rate & 0.5 pep/(nL $\cdot$ day)\smallskip\\
    $\delta_A$ & antigen presentation degradation rate & 0.8 day$^{-1}$\smallskip\\ % 4.38

    % INFLAMMATORY SIGNALS
    $r_I$ & inflammation source rate & 0.4 (ng$\cdot$nL)/(cells$^2$ $\cdot$ day)\smallskip\\
    $\delta_I$ & inflammation degradation rate & 3.0 day$^{-1}$\smallskip\\     % 0.01
    
    \hline\noalign{\smallskip}
    % T cells
    $r_E$ & effector T cell growth coefficient & 1.0 day$^{-1}$\smallskip\\
    $E^*$ & effector T cell base steady state & 5.0 cells/nL\smallskip\\
    $r_S$ & non-effector T cell growth coefficient & 1.0 day$^{-1}$\smallskip\\
    $S^*$ & non-effector T cell base steady state & 5.0 cells/nL\smallskip\\
    $\beta$ & effector T cell recruitment coefficient & 0.009 nL$^3$/(pep$\cdot$ng$\cdot$cells$\cdot$day)\smallskip\\
    $\gamma$ & non-effector T cell recruitment coefficient & 37.414 nL/(cells $\cdot$ day)\smallskip\\
    \hline
  \end{tabular}
  \caption{Summary of \bb{baseline parameter} estimates.}\label{table:baselineparams}
\end{table}

\begin{table}[htbp]
\centering
    \begin{tabular}{l l l l l}
    \hline\noalign{\smallskip} 
    Variable & Initial value\\
    \noalign{\smallskip}\thickhline\noalign{\smallskip}
    $C$ & steady state $C^*$ \smallskip\\
    $A$ & 0 peptides/nL \smallskip\\
    $I$ & 0 ng/nL \smallskip\\
    $E$ & 0 cells/nL \smallskip\\
    $S$ & base steady state $S^*$ \smallskip\\
    \hline
  \end{tabular}\smallskip
  \caption{\bb{Summary of initial values}.}\label{table:ICs}
\end{table}

%%%%%%%%%%%%%%%%%%%%%%%%%%%%%%%%%%%%%%%%%%%%%%%%%%%%%%%%%%%%%%%%%%%%%
%%%%%%%%%%%%%%%%%%%%%%%%%%%%%%%%%%%%%%%%%%%%%%%%%%%%%%%%%%%%%%%%%%%%%

%%% CANCER PARAMS %%%f
\noindent \ii{Parameters for cancer, $C$}:\\[1mm]
Carlson (2003) estimated the doubling time of metastatic tumours to be 8 to 212 days~\cite{carlsonTumorDoublingTime2003}. We take the maximum relative growth rate $r_{\text{max}} = (\ln 2)/(8 \text{days}) \approx 0.09$/day and fit the parameter $r_C \approx 30$ with MATLAB. Next, we estimate the steady state cancer concentration, $C^*$. The volume of a 1 cm$^3$ tumour is typically assumed to have between 100 million to 1 billion cancer cells~\cite{delmonteDoesCellNumber2009, narodDisappearingBreastCancers2012}. For simplicity, we assume the size of the tumour microenvironment to be the size of the tumour and set $C^* = 10^9\,\text{cells}/1\,\text{cm}^3 = 10^3\,\text{cells}/\text{nL}$. Finally, we estimate the mass-action killing rate of cancer cells, $\kappa$, by effector T cells. We assume that the killing rate is a proportion of the mass-action rate at which effector T cells contact cancer cells, denoted $r_\text{con1}$. This mass-action contact rate is calculated using the actual contact rate between the T cells and cancer cells, denoted $r_\text{con2}$, as follows. By the law of mass action, we can write

\begin{align*}
    \text{(mass-action contact rate)}[\text{effector T}][\text{cancer}] &= \text{(contact rate)}[\text{cancer}]\\[2mm]
    r_\text{con1}\,K\,C &= r_\text{con2}\,C.
\end{align*}

Rearranging, the mass action contact rate $r_\text{con1} = r_\text{con2}/K$. To calculate the right hand side, we use the estimate that in the lymph node of Catron et al. (2014), one T cell and one DC will have 0.20 $\pm$ 0.06 interactions per hour, or 4.8 $\pm$ 1.4 interactions per day~\cite{catronVisualizingFirst502004}. Moreover, the T cell concentration in this scenario is 500 cells per 4.2 $\mu$L-sized lymph node, or 119 cells/$\mu$L. This implies the mass-action contact rate for Catron et al. is 4.8/119 = 0.04 $\mu$L/(cells $\cdot$ day) = 40 nL/(cells $\cdot$ day). For simplicity, we assume that the T cell to DC mass-action interaction rate in the lymph to be comparable to the T cell to cancer cell mass-action interaction rate in the tumour microenvironment and take $r_\text{con1}$ to be 40 nL/(cells $\cdot$ day).
It is unlikely that every antigen-specific T cell-cancer interaction leads to a successful apoptosis of the cancer. In particular, cancer leverages immunosuppressive mechanisms that can suppress the cytotoxicity of CD4$^+$ T cells mobilised to destroy them~\cite{medemaBlockadeGranzymePerforin2001} and we assume that T cell killing might only be successful after repeated interactions. In our model, we assume the tumour is immune evasive and that the probability of a successful kill to be 3\% per interaction. We take the mass-action kill rate $\kappa = 0.03 \times r_\text{con1} = 0.03\times 40 = 1.2$ nL/(cells $\cdot$ day). The choice of 3\% killing probability was chosen to give interesting dynamics $\kappa$'s effect on the results will be discussed in a sensitivity analysis.\\

%%% ANTIGEN PRESENTATION and INFLAMMATION PARAMETERS %%%%
\noindent \ii{Parameters for the signals $A$ and $I$}:\\[1mm]
Belz et al. (2007) shows that the level of antigen presentation following the third day of an infection decays with a half-life between 19 h and 20.4 h~\cite{belzKillerCellsRegulate2007}. Using a half-life of 20 h, we set the antigen presentation degradation rate $\delta_A = (\ln 2)/20$ h$^{-1} \approx$ 0.8 day$^{-1}$. The other parameters, $r_A = 0.5$, $r_I = 0.4$, $\delta_I = 3.0$ were fitted for simulations.\\

%%% EFFECTOR AND NON-EFFECTOR PARAMETERS %%%%
\noindent \ii{Parameters for T cells, $E$ and $S$}:\\[1mm]
The parameters $r_E = r_S = 1$, $E^* = S^* = 5$, $\beta \in [0.008, 0.009]$ and $\gamma \in [37.414, 37.5]$ were fitted for simulations.\\[1mm]

%%% ICs %%
\noindent \ii{Initial values}:\\[1mm]
For simplicity, we assume the cancer has reached steady state and set $C(0) = C^*$. We also assume that tumour escape has occurred so that $E(0) \approx 0$ and $S(0) \approx S^*$, which broadly correspond to an immunosuppressive environment where tumour infiltrating lymphocytes (TILs) have largely been suppressed~\cite{antoheTumorInfiltratingLymphocytes2019}. Since the variables for antigen presentation and inflammation, $A$ and $I$, do not reveal obvious steady states, for simplicity, we set their initial levels at some low level $A(0) = 1$ and $I(0) = 1$.

%%%%%%%%%%%%%%%%%%%%%%%%%%%%%%%%%%%%%%%%%%%%%%%%%%%%%%%%%%%%%%%
%%%%%%%%%%%%%%%%%%%%%%%%%%%%%%%%%%%%%%%%%%%%%%%%%%%%%%%%%%%%%%%
%%%%%%%%%%%%%%%%%%%%%%%%%%%%%%%%%%%%%%%%%%%%%%%%%%%%%%%%%%%%%%%%
%%%%%%%%%%%%%%%%%%%%%%%%%%%%%%%%%%%%%%%%%%%%%%%%%%%%%%%%%%%%%%%
%%%%%%%%%%%%%%%%%%%%%%%%%%%%%%%%%%%%%%%%%%%%%%%%%%%%%%%%%%%%%%%
%%%%%%%%% 4. RESULTS %%%%%%%%%%%%%%%%%%%%%%%%%%%%%%%%%%%%%%%%%%
%%%%%%%%%%%%%%%%%%%%%%%%%%%%%%%%%%%%%%%%%%%%%%%%%%%%%%%%%%%%%%%
%%%%%%%%%%%%%%%%%%%%%%%%%%%%%%%%%%%%%%%%%%%%%%%%%%%%%%%%%%%%%%%%
%%%%%%%%%%%%%%%%%%%%%%%%%%%%%%%%%%%%%%%%%%%%%%%%%%%%%%%%%%%%%%%
%%%%%%%%%%%%%%%%%%%%%%%%%%%%%%%%%%%%%%%%%%%%%%%%%%%%%%%%%%%%%%%
%%%%%%%%%%%%%%%%%%%%%%%%%%%%%%%%%%%%%%%%%%%%%%%%%%%%%%%%%%%%%%%

\newpage
\section{Results}\label{sec:results}
\begin{figure}
\centering
\includegraphics[scale=1]{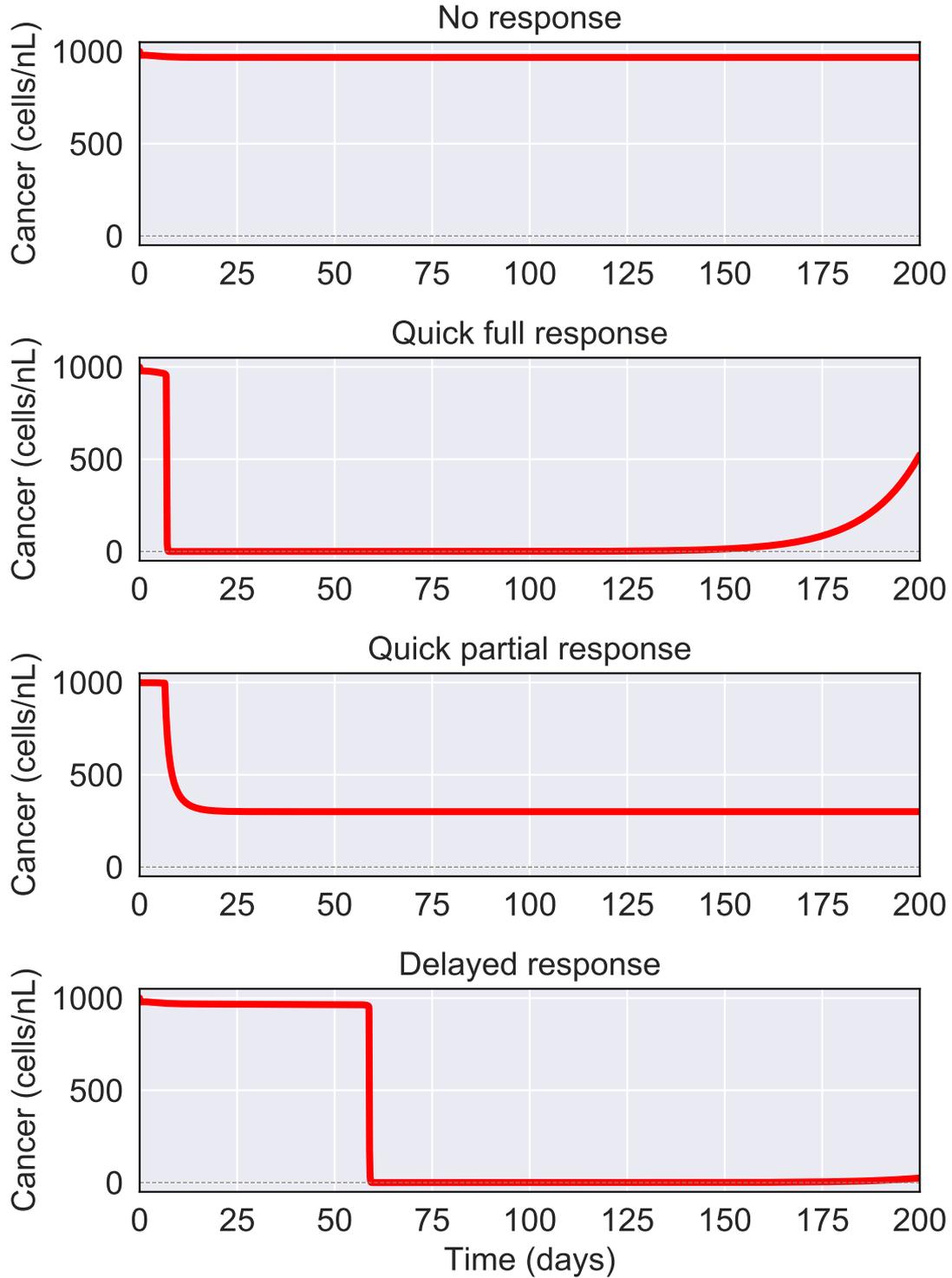}
\caption{\bb{Qualitative responses generated by the model}. The four cases are: 1) an absence of response, with very limited anti-cancer immune activity; 2) a quick and full response, where the tumour is almost eradicated in the span of of weeks; 3) a quick and partial response, where the tumour volume partially recedes; and 4) delayed responses, where anti-cancer activity suddenly escalates after months of minimal activity. All response types can be generated by four parameters: $\beta$, $\gamma$, $E^*$ and $r_\text{max}$. See Table \ref{table:baselineparamsallcases} for these parameter sets. Note that all simulations are idealised in that emphasis is placed on broadly generating the different response cases and contrasting their differences, since datasets for the delayed-response case are generally limited.}
\label{SpikeODE_result_responses_all}
\end{figure}

%%% summary
Our model qualitatively recreates several experimentally observed outcomes seen with immune checkpoint inhibitors. In Fig.~\ref{SpikeODE_result_responses_all}, we collate the idealised simulations seen across the four response cases generated by our model. These cases are as follows: 1) an \bb{absence of response} with very limited anti-cancer immune activity; 2) a \bb{quick and full response}, where the tumour is seemingly eradicated in the span of weeks; 3) a \bb{quick and partial response}, where the tumour volume partially recedes; and 4) a \bb{delayed responses} where anti-cancer activity suddenly escalates after months of minimal activity~\cite{wolchokGuidelinesEvaluationImmune2009}. To maintain simplicity for comparing and contrasting response cases, we note our model is not fit to any particular dataset, which are scarce and perhaps not representative.\\

\begin{figure}
\centering
\includegraphics[scale=1]{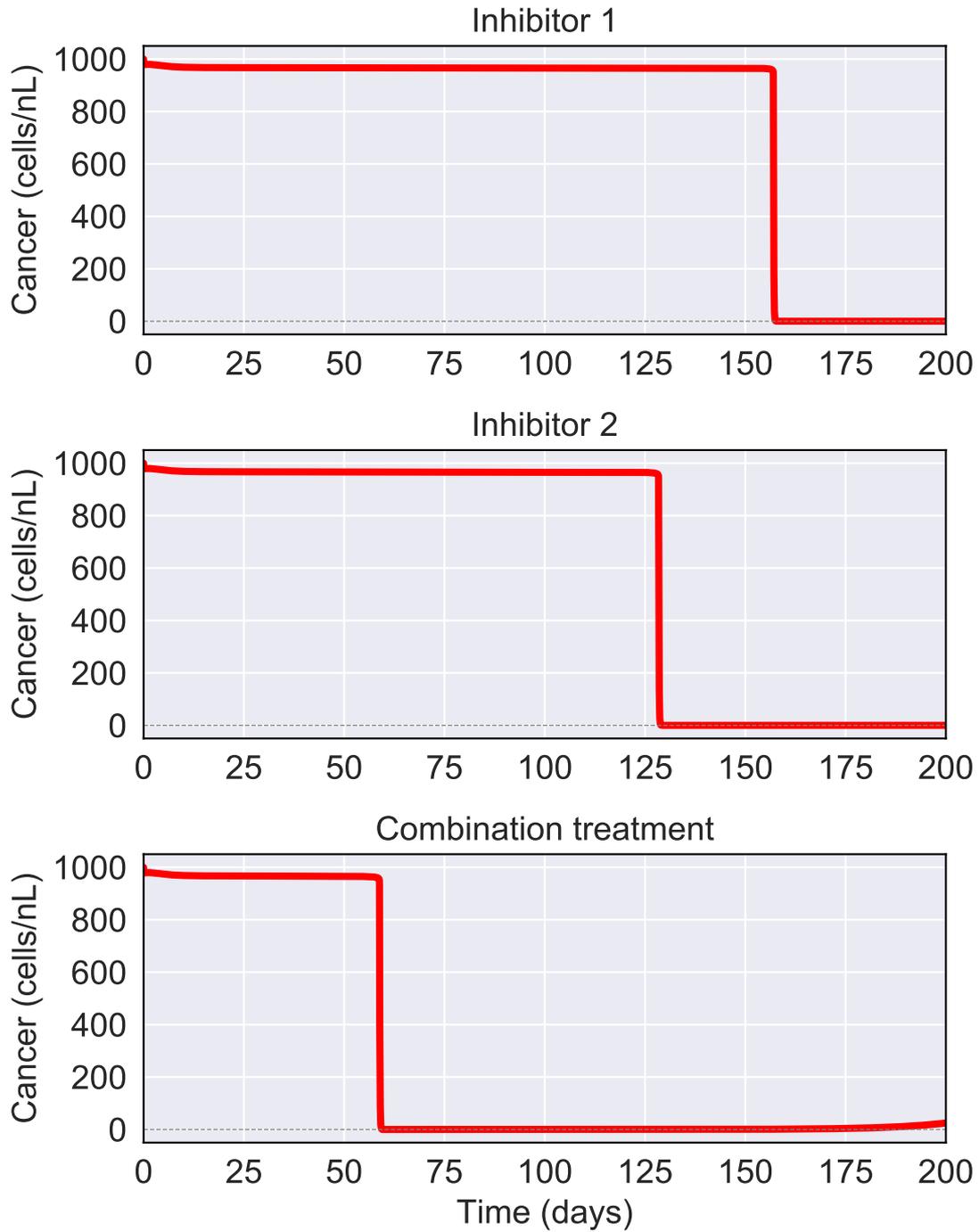}
\caption{\bb{Delayed responses}. One of the responses produced by our model is the delayed-response case, where the tumour volume does not react for a period spanning months or years. This delay can be calibrated continuously by perturbing $\beta$ and $\gamma$, which correspond to the level of CTLA-4 and PD-1 respectively. The top graph represents the administration of only a CTLA-4 inhibitor ($\beta = 0.009$) and the middle graph represents the administration of only a PD-1 inhibitor ($\gamma = 37.414$). These individual inhibitors generate delay lengths of 5 and 4 months respectively. The administration of both inhibitors ($\beta = 0.009, \gamma = 37.414$) results in a synergistic benefit to the delay length, reducing it down to 2 months. The precise parameter values can be found in Table \ref{table:baselineparamsallcases}.}
\label{SpikeODE_result_responses_delayed_only}
\end{figure}

\begin{figure}
\centering
\includegraphics[scale=1]{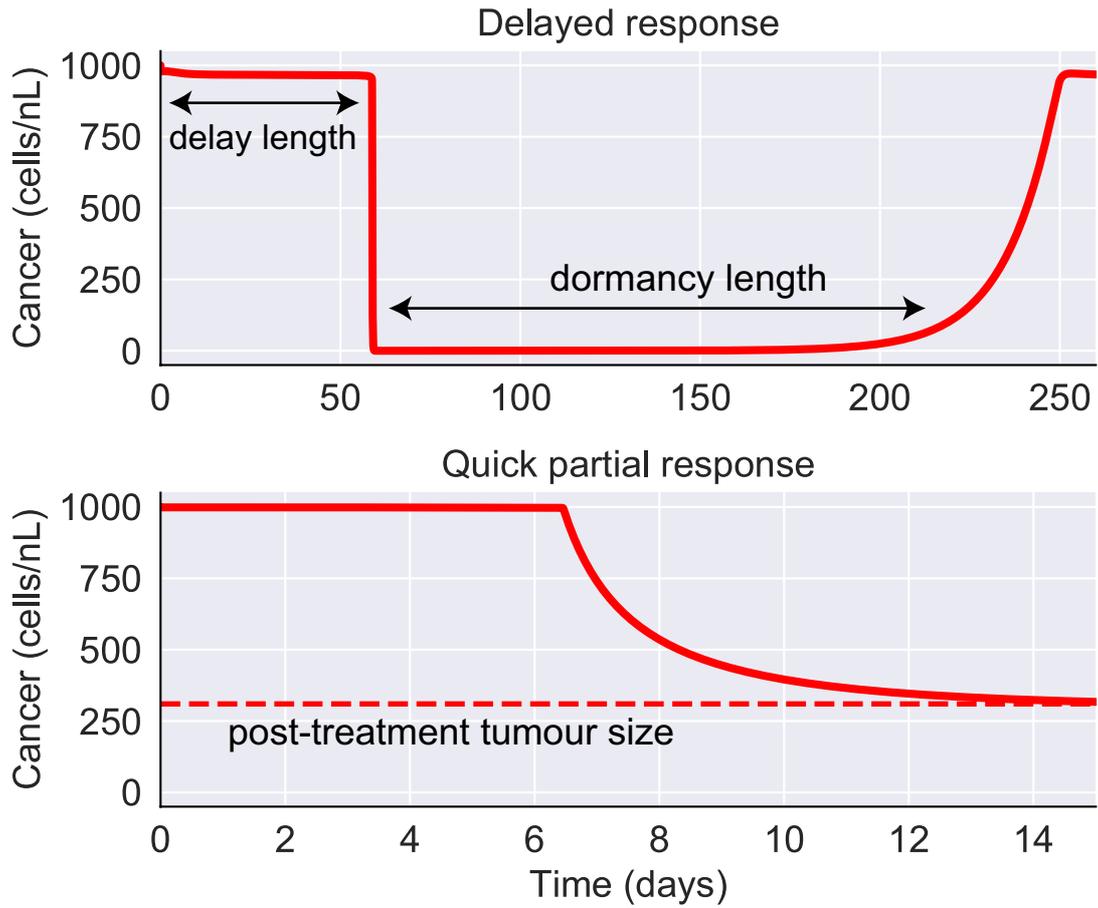}
\caption{\bb{Quantities of interest}. Our model predicts that for delayed responses, the cancer eventually returns after a period of dormancy. Two quantities of interest are, therefore, the \bb{delay length} and \bb{dormancy length}. Moreover, our model predicts that in the case of a partial response, where the tumour size stabilises, the delay is always quick. Thus for these quick partial responses, another insightful quantity is the \bb{post-treatment tumour size}.}
\label{SpikeODE_result_clinical_quantities}
\end{figure}

\begin{figure}
\centering
\includegraphics[scale=1]{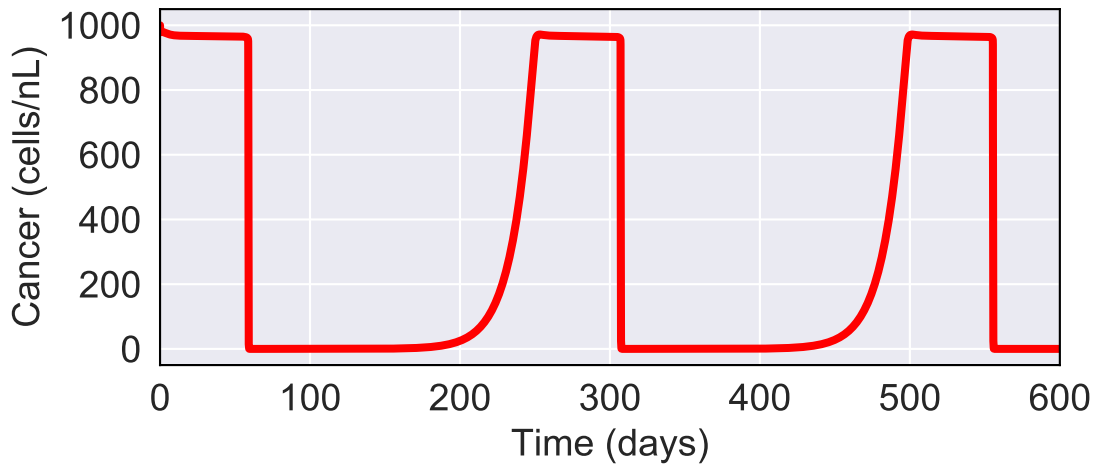}
\caption{\bb{Cyclic behaviour of delayed responses}. For the delayed response case, the cancer burden is reduced to near zero, but not completely eradicated. This later results in relapse. This behaviour repeats indefinitely.}
\label{SpikeODE_result_cyclic}
\end{figure}

\begin{figure}
\centering
\includegraphics[scale=1]{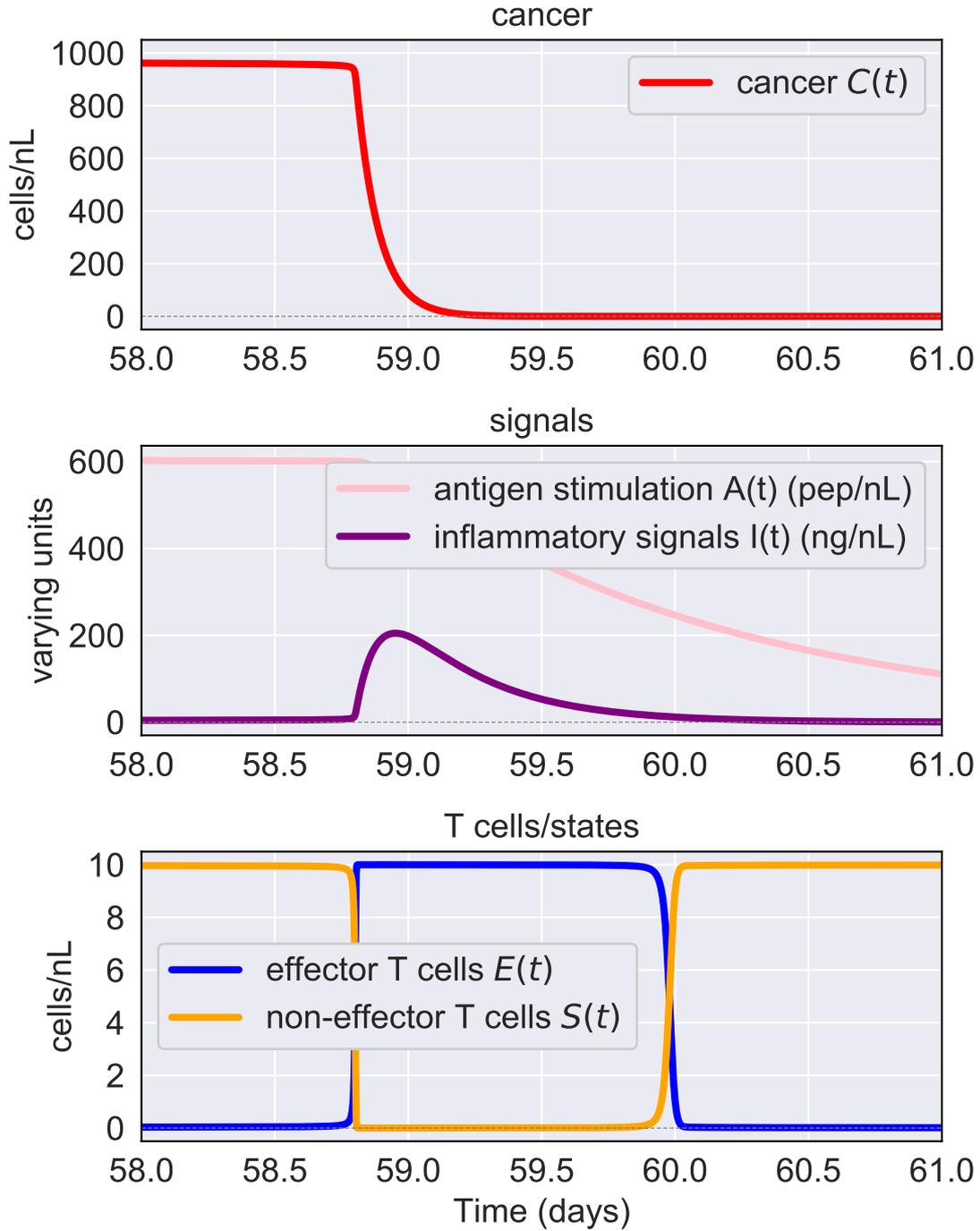}
\caption{\bb{T-cell response dynamics}. This figure plots all five variables in the ODE: $C$, $A$, $I$, $E$ and $S$ for the delayed-response case around the time the delay occurs. Here, we see that the levels of effector and non-effector T cells, $E$ and $S$, sea-saw in the tumour microenvironment. However, the effector cells are short-lived as depletion of antigen and inflammatory signals caused by eradication of the tumour quickly shuts off recruitment for new $E$ cells. Note that even after the $E$ cells become suppressed, the tumour stays suppressed for a longer period of time. In Fig. \ref{SpikeODE_result_T_cell_responses_long}, one can see that the initial $E$ response is capable of suppressing the tumour for a very long period of time. This set of simulations is generated from the baseline parameters that simulate a 2-month delay---see Table \ref{table:baselineparams}.}
\label{SpikeODE_result_T_cell_responses_zoomed}
\end{figure}

\begin{figure}
\centering
\includegraphics[scale=1]{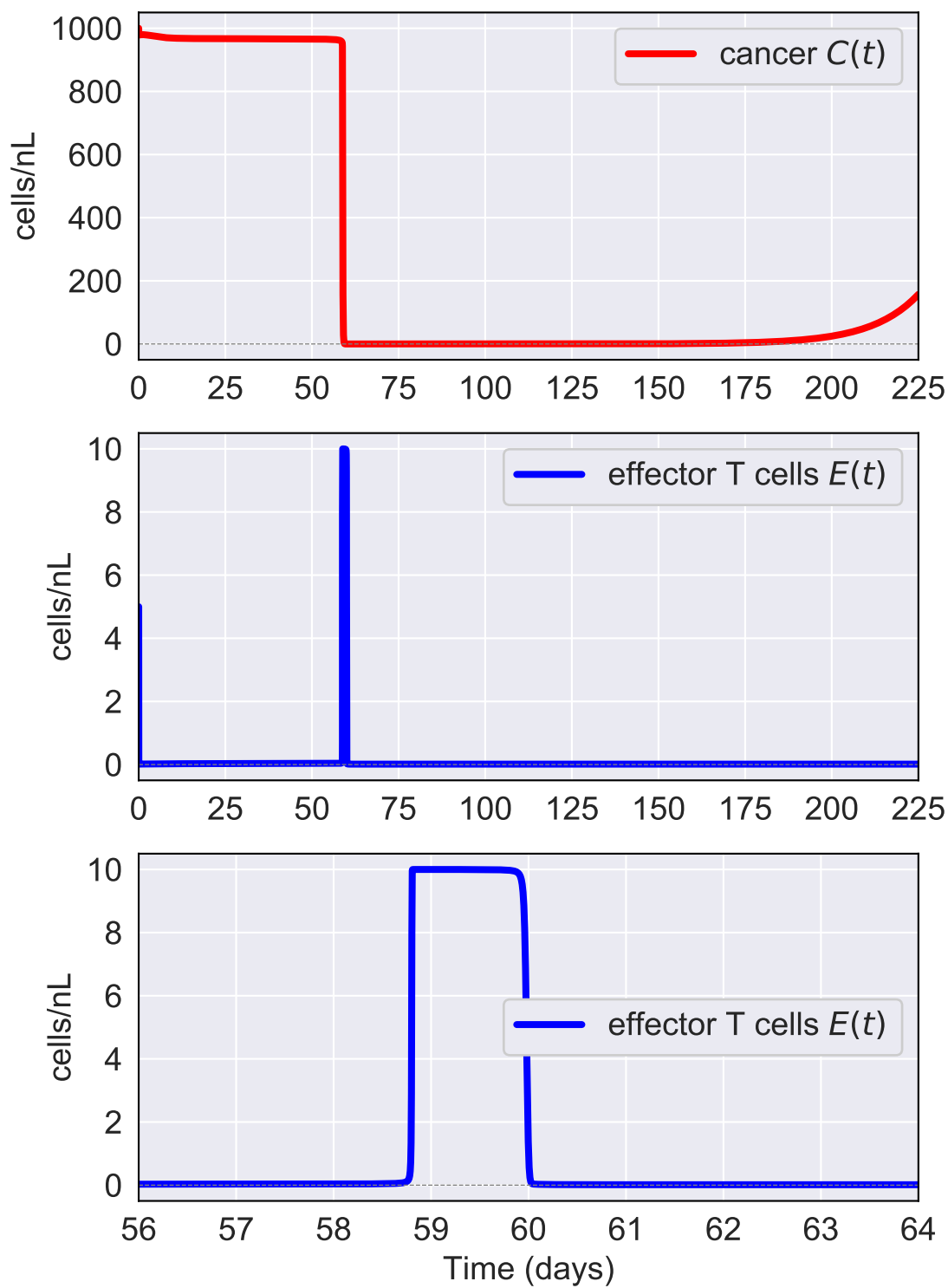}
\caption{\bb{Effector T cell response and tumour suppression}. Our model suggests that a relatively short effector T cell response (spanning days) can suppress the tumour for a much longer period (months). The middle and zoomed-in bottom graph shows an anti-cancer T cell response spanning less than two days is able to hold the cancer at bay for roughly 6 months (top graph). This simulation is based on the baseline parameters that simulate a 2-month delay length---see Table \ref{table:baselineparams}.}
\label{SpikeODE_result_T_cell_responses_long}
\end{figure}

For the delayed response case, the model is able to replicate delays of varying lengths, spanning months or even years. Specifically, we tune the parameters $\beta$ and $\gamma$, linked to the checkpoint immune checkpoint receptors CTLA-4 and PD-1, to give a delayed response after 5 and 4 months, respectively. Moreover, our model replicates the synergistic benefit seen by combining anti-CTLA-4 (Ipilimumab) and PD-1 (Nivolumab) inhibitors~\cite{khairCombiningImmuneCheckpoint2019}. Plots of these delayed responses can be found in Fig.~\ref{SpikeODE_result_responses_delayed_only}. Given its clinical significance, we shall henceforth refer to the length of time between commencement of treatment and the delayed anti-tumour T cell response as the \ii{delay length} (see Fig. \ref{SpikeODE_result_clinical_quantities}). The parameter sets for these different responses are discussed in the parameter sensitivity section.\\
\\
\indent The long-term behaviour of the model falls into two categories: \bb{tumour equilibrium} and \bb{tumour dormancy}. Tumour equilibrium occurs during the partial quick response case. Here, the anti-cancer T cell response suppresses the tumour volume to some reduced steady state, which holds indefinitely. Tumour dormancy occurs during the delayed response case. Here, the tumour volume is suppressed close to zero, which does not hold indefinitely and the tumour later relapses. This behaviour then repeats cyclically (see Fig. \ref{SpikeODE_result_cyclic}). In general, our model exhibits cyclic behaviour in all response cases except for the quick partial response, where the cancer responds quickly but stabilises away from the cancer-free equilibrium--instead stabilising somewhere between $C \approx 0$ and $C = C^*$.\\
\\
\indent Alongside \ii{delay length}, two additional clinically relevant quantities that arise from considering the long-term behaviour of the model are \ii{post-treatment tumour size} and \ii{dormancy length}. A summary of these three quantities can be found in Fig.~\ref{SpikeODE_result_clinical_quantities}. We note that in all cases---response or no response, delayed or not---our model does not predict the permanent eradication of the tumour.\\
\\
\indent We observe that the two T cell populations, or \ii{states}, $E$ and $S$ compete in a zero-sum game in the tumour microenvironment. Specifically, an elevation in effector T cell levels, $E$, is accompanied by a reduction in the level for $S$, and vice versa. This sea-saw behaviour between $E$ and $S$ can be observed most clearly in the delayed response simulations (see Fig.~\ref{SpikeODE_result_T_cell_responses_zoomed}), where a strong anti-cancer effector T cell response occurs in tandem with a similarly potent suppression of non-effector T cells at the tumour site.\\ 
\\
\indent This interplay between $E$ and $S$ draws loose parallels to the interactions between effector T cells (Teff) and regulatory T cells (Treg) within tumours. Here, experiments have shown immune checkpoint inhibitors to be effective in increasing the ratio of Teff/Treg, resulting in rapid tumour rejection \cite{oleinikaSuppressionSubversionEscape2013, quezadaCTLA4BlockadeGMCSF2006, curranTumorVaccinesExpressing2009}. Moreover, our simulations suggest that the duration of time T cells spent in an effector state (i.e. as $E$ cells) is comparatively short. We suggest the model offers a compelling explanation: because the cancer eradication is swift, the antigen and inflammation levels required to sustain continued $E$ recruitment also decreases rapidly (see Fig.~\ref{SpikeODE_result_T_cell_responses_zoomed}). However, it is less obvious why the dormancy length is significantly longer than the length of time T cells stay in effector state (see Fig.~\ref{SpikeODE_result_T_cell_responses_long}). This suggests although the T cell response is sharp and short, it may be sufficient to prevent relapse for a much longer period of time.\\
\\
\indent Finally, we note from Fig. \ref{SpikeODE_result_T_cell_responses_zoomed} that the total T cell population appears to be conserved, where $E(t) + S(t) = E^* + S^*$ after an initial transient. This can be explained as follows. Assuming $r_E = r_S = 1$ and adding Equations \ref{eq:spikeE} and \ref{eq:spikeS} together gives
    \begin{align*}
        \frac{dE}{dt} + \frac{dS}{dt} &= k(E + S - (E^* + S^*)).
    \end{align*}
    This is a linear differential equation in $E + S$ with the solution $E + S = E^* + S^* + A\,e^{-kt}$ for some arbitrary constant $A$. Our initial conditions $E(0) = 0$ and $S = S^*$ render $A = 0$. For other initial conditions, the exponential decays quickly on a timescale in months, resulting in $E + S \approx E^* + S^*$ after an initial transient.

% %%%%%%%%%%%%%%%%%%%%%%%%%%%%%%%%%%%%%%%%%%%%%%%%%%%%%%%%%%%%%%%%%%%%%%%%%%%%
% Res DIAGRAM: Responses all
% %%%%%%%%%%%%%%%%%%%%%%%%%%%%%%%%%%%%%%%%%%%%%%%%%%%%%%%%%%%%%%%%%%%%%%%%%%%%

% %%%%%%%%%%%%%%%%%%%%%%%%%%%%%%%%%%%%%%%%%%%%%%%%%%%%%%%%%%%%%%%%%%%%%%%%%%%%%
% Res DIAGRAM: Responses delayed only
% %%%%%%%%%%%%%%%%%%%%%%%%%%%%%%%%%%%%%%%%%%%%%%%%%%%%%%%%%%%%%%%%%%%%%%%%%%%%

% %%%%%%%%%%%%%%%%%%%%%%%%%%%%%%%%%%%%%%%%%%%%%%%%%%%%%%%%%%%%%%%%%%%%%%%%%%%%
% Res DIAGRAM: Clinically-relevant quantities of delay length and dormancy length
% %%%%%%%%%%%%%%%%%%%%%%%%%%%%%%%%%%%%%%%%%%%%%%%%%%%%%%%%%%%%%%%%%%%%%%%%%%%%

% %%%%%%%%%%%%%%%%%%%%%%%%%%%%%%%%%%%%%%%%%%%%%%%%%%%%%%%%%%%%%%%%%%%%%%%%%%%%%
% Res DIAGRAM: Delayed cyclic
% %%%%%%%%%%%%%%%%%%%%%%%%%%%%%%%%%%%%%%%%%%%%%%%%%%%%%%%%%%%%%%%%%%%%%%%%%%%%

% %%%%%%%%%%%%%%%%%%%%%%%%%%%%%%%%%%%%%%%%%%%%%%%%%%%%%%%%%%%%%%%%%%%%%%%%%%%%
% Res DIAGRAM: T cell responses zoomed in
% %%%%%%%%%%%%%%%%%%%%%%%%%%%%%%%%%%%%%%%%%%%%%%%%%%%%%%%%%%%%%%%%%%%%%%%%%%%%

% %%%%%%%%%%%%%%%%%%%%%%%%%%%%%%%%%%%%%%%%%%%%%%%%%%%%%%%%%%%%%%%%%%%%%%%%%%%%
% Res DIAGRAM: T cell responses long
% %%%%%%%%%%%%%%%%%%%%%%%%%%%%%%%%%%%%%%%%%%%%%%%%%%%%%%%%%%%%%%%%%%%%%%%%%%%%

%%%%%%%%%%%%%%%%%%%%%%%%%%%%%%%%%%%%%%%%%%%%%%%%%%%%%%%%%%%%%%
%%%%%%%%%%%%%%%%%%%%%%%%%%%%%%%%%%%%%%%%%%%%%%%%%%%%%%%%%%%%%%
%%%%%%%%%%%%%%%%%%%%%%%%%%%%%%%%%%%%%%%%%%%%%%%%%%%%%%%%%%%%%%
%%%%%%%%%%%%%%%%%%%%%%%%%%%%%%%%%%%%%%%%%%%%%%%%%%%%%%%%%%%%%%
%%%%%%%%%%%%%%%%%%%%%%%%%%%%%%%%%%%%%%%%%%%%%%%%%%%%%%%%%%%%%%
%%%%%%%%% 5. SENSITIVITY ANALYSIS %%%%%%%%%%%%%%%%%%%%%%%%%%%%
%%%%%%%%%%%%%%%%%%%%%%%%%%%%%%%%%%%%%%%%%%%%%%%%%%%%%%%%%%%%%%%
%%%%%%%%%%%%%%%%%%%%%%%%%%%%%%%%%%%%%%%%%%%%%%%%%%%%%%%%%%%%%%%
%%%%%%%%%%%%%%%%%%%%%%%%%%%%%%%%%%%%%%%%%%%%%%%%%%%%%%%%%%%%%%%
%%%%%%%%%%%%%%%%%%%%%%%%%%%%%%%%%%%%%%%%%%%%%%%%%%%%%%%%%%%%%%%
%%%%%%%%%%%%%%%%%%%%%%%%%%%%%%%%%%%%%%%%%%%%%%%%%%%%%%%%%%%%%%%

\newpage
\section{Sensitivity of parameters}\label{sec:sensitivity}
%%%
Recall that our model generates four response outcomes: no response, quick full response, quick partial response and delayed response (see Fig. \ref{SpikeODE_result_responses_all}). The majority of parameter sets generate no response or quick full responses. In contrast, generating a delayed outcome is difficult and requires careful selection of parameter values. Upon obtaining a finely calibrated set of parameters that gives a delayed response, the delay can easily vanish (giving a quick full response) or extend to infinity (resulting in no response) with small perturbations of the parameters. Moreover, similarly to delayed responses, the quick partial response is not easily obtained, although it can be derived from the delayed response by elevating the cancer growth rate to a certain level $r_{\text{max}} = 1$. Since this exceeds biological feasibility, we regard the full partial response case as unlikely, and justify its inclusion into the results section on the basis of completeness. 

\vv{0.5pc}

In summary, all four responses can be generated from delayed-response parameters as a starting base. In particular, a baseline parameter set that generates a two-month delay intended to simulate a combination anti-CTLA-4 (Ipilimumab) and PD-1 (Nivolumab) can be found in Table \ref{table:baselineparams}. From this baseline, we can generate delays of varying durations (months to years) alongside the different types of responses (no response to quick response) by perturbing four parameters: $\beta$, which governs the effect of CTLA-4 and its inhibitor; $\gamma$, which governs the effect of PD-1 and its inhibitor;  $E^*$, which calibrates a patient's individual T cell response to cancer, such as their level of tumour-infiltrating lymphocytes~\cite{hendryAssessingTumorInfiltrating2017, yiSynergisticEffectImmune2019}; and $r_\text{max}$, which calibrates cancer growth rate. A summary of these parameters can be found in Table~\ref{table:baselineparamsallcases}.

\vv{0.5pc} 

\begin{figure}
\centering
\includegraphics[scale=0.8]{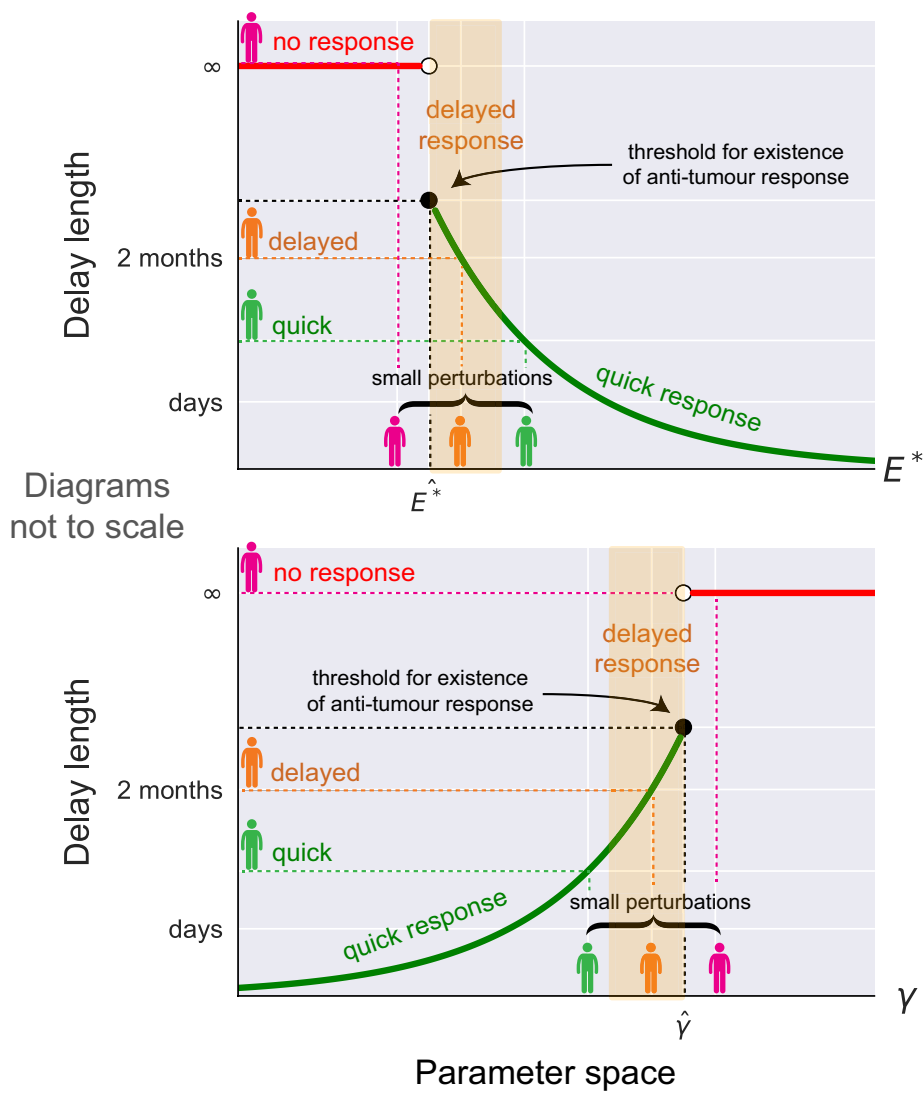}
\caption{\bb{Illustrative diagram of response time versus model parameters}. \ii{Note: diagrams are not to scale.} In our model, the relationship between the length of the delay before an anti-tumour response and most parameters is non-linear and discontinuous. Take, for instance, response time versus the effector T cell parameter $E^*$ (top graph). For values of $E^* < \hat{E^*}$, the model predicts an absence of anti-tumour response. This could be due to an insufficient level of tumour-infiltrating lymphocytes in the patient's tumour. For $E^* > \hat{E^*}$, the non-linear curve (not to scale) gives rise to a thin region for $E^*$ that results in a delayed response spanning months. Beyond this, the response time is quick. Note that for our baseline parameters (see Table \ref{table:baselineparams}), $\hat{E^*} = 5.0$. As another example, take the relationship between the response time and the PD-1 associated parameter $\gamma$. As $\gamma$ increases, the model predicts that immunosuppressive pressures on effector T cells increase at an expanding rate. This results in a thin region of $\gamma < \hat{\gamma}$ that generates a delayed response spanning months. At the threshold $\gamma = \hat{\gamma}$, the delay vanishes and the response time becomes infinitely long. Note that for our baseline parameters, $\hat{\gamma} \approx 37.42$.
\bb{Patient variations and immune profiles}. Broadly speaking, parameter combinations can be considered as mapping to immune profiles belonging to different people. For instance, perturbations of $E^*$ might result in slightly different levels of TILs in a tumour between two individuals, while perturbations of $\gamma$ might result in different levels of PD-1. These parameter-delay plots shine light into why two different patients may exhibit wildly different response times. In particular, two patients with immune profiles that position them on different sides of the jump discontinuity will affect how they react to the tumour and immunotherapy. Figure~\ref{SpikeODE_result_delay_vs_parameters_drugs} describes the various treatment scenarios possible with our model.}

\label{SpikeODE_result_delay_vs_parameters_v4}
\end{figure}

For delayed responses, we measured the sensitivity of the two clinically relevant quantities, \ii{delay length} and \ii{dormancy duration}, to our model parameters, by perturbing parameters one at a time while holding the others at baseline. The results can be found in Table~\ref{table:parameterOATsensitivity}. We found that the dormancy length is sensitive only to the maximum cancer growth rate, $r_{\text{max}}$, while the delay length is highly sensitive to all other parameters. That is, starting from our baseline parameters that generate a 2-month delay, almost every single parameter has the power to substantially affect the length of this delay. For instance, a 1\% change in our CTLA-4 associated parameter $\beta$---reflecting a change from one patient to another---can shorten the delay length by 75\%, putting their response time to treatment into days rather than months. See Fig. \ref{SpikeODE_result_delay_vs_parameters_v4} for an illustration.

\vv{0.5pc}

In regard to which parameters lengthen or shorten the delay, we have that the cancer killing rate by T cells $\kappa$, signal degradation rates $\delta_A$ and $\delta_I$, non-effector T cell steady state $S^*$ and recruitment rate $\gamma$ (associated with PD-1) induce a longer delay. None of these relationships are surprising with the exception of the killing rate: we unexpectedly observe that higher cytotoxic T cells correlate with a longer delay prior to a robust anti-cancer response, which highlights the complex dynamics underlying tumour-immune interactions. In contrast, the following parameters correlate with shorter delay lengths: cancer growth rates $r_C$ and $r_{\text{max}}$, cancer size $C^*$, positive signal source rates $r_A$ and $r_I$, intrinsic T cell growth rates $r_E$ and $r_S$, effector T cell steady state $E^*$ and recruitment rate $\beta$ (associated with CTLA-4). Some of these relationships perhaps also run counter to expectations; for instance the notion that larger and faster-growing tumours (higher $C^*$ and $r_C$) coincide with a speedier treatment response. Therefore our model suggests dynamic and at-times surprising relationships between treatment delay length and the parameters.

\vv{0.5pc}

Moreover, the relationship between the delay length and general parameter is non-linear and discontinuous in a way that illustrates how unconventional a delayed response might be. For instance, consider the PD-1 and immunosuppression-associated parameter $\gamma$. We observe that the delay length increases at an expanding rate as $\gamma$ increases linearly. As $\gamma$ approaches a critical value of approximately $\hat{\gamma} \approx 37.42$, the delay length increases to a very fast rate, ultimately jumping from a finite amount of time to becoming infinitely long. In particular, the model admits only a thin region in the parameter space near $\hat{\gamma}$ that generates a delayed response spanning into the months, implying that only patients with specific immune system profiles may experience delayed responses to treatment. This non-linear, discontinuous relationship can be seen in Fig. \ref{SpikeODE_result_delay_vs_parameters_v4}.

\vv{0.5pc}

Finally, having seen the simulations, we are in a better position to comment on the choice for the cancer growth rate, $f(C)$. In particular, for the delayed-response case, the saturating function $f(C)$ satisfies several properties: it permits logistic-like growth with $f(C) := r_C\,(1 - C/C^*)$ during the delay phase, whilst ensuring the cancer never grows too quickly with $f(C) := r_\text{max}$ during the cancer killing phase. Without the latter condition, the cancer can grow at a biologically unreasonable rate due to our parameter fits, notably the rate parameter $r_C = 30$. Other saturating growth rates $f(C)$ are possible; however, our choice of $f(C)$ was effective in balancing a set of requirements.

% %%%%%%%%%%%%%%%%%%%%%%%%%%%%%%%%%%%%%%%%%%%%%%%%%%%%%%%%%%%%%%%%%%%%%%%%%%%%%
% Sens TABLE: parameter values for 4 types of responses and 3 delayed lengths. BASELINE PERTURBATION FOR ALL CASES
% %%%%%%%%%%%%%%%%%%%%%%%%%%%%%%%%%%%%%%%%%%%%%%%%%%%%%%%%%%%%%%%%%%%%%%%%%%%%%

\begin{table}[htbp]
  \centering
  \begin{tabular}{l l l l l}
    \hline\noalign{\smallskip} 
    Types of responses & $\beta$ & $\gamma$ & $E^*$ & $r_\text{max}$\\
    \noalign{\smallskip}\thickhline\noalign{\smallskip}
    No response & $0.0089988$ & $37.4168$ & 5 & 0.09 \smallskip\\
    Quick full response & $0.009$ & $37.4168$ & 5.5 & 0.09 \smallskip\\
    Quick partial response & $0.0089988$ & $37.414$ & 5 & 1\smallskip\\
    Delayed response & $0.009$ & $37.414$ & 5 & 0.09 \smallskip\\
    \hline
  \end{tabular}

\vv{0.5pc}
  
  \hh{-1pc}
  \begin{tabular}{l l l l l l}
    \hline\noalign{\smallskip} 
    Types of delayed responses & $\beta$ & $\gamma$\\
    \noalign{\smallskip}\thickhline\noalign{\smallskip}
    No treatment & $0.0089988$ & $37.4168$ \smallskip\\
    Inhibitor 1 (5 months) & $0.009$ & $37.4168$ \smallskip\\
    Inhibitor 2 (4 months) & $0.0089988$ & $37.414$ \smallskip\\
    \bb{Combination (2 months)**} & $0.009$ & $37.414$ \smallskip\\
    \hline
  \end{tabular}\smallskip
  \caption{Different response types and the key parameters that generate them. Only four parameters ($\beta$, $\gamma$, $E^*$ and $r_\text{max}$) are needed to generate all different response types (quick full, quick partial, delayed and no response). For the delayed responses, only $\beta$ and $\gamma$ (associated with CTLA-4 and PD-1) are needed. Our model predicts that individual inhibitors work synergistically in combination to result in a quicker response. This has been widely observed in clinical trials.
 **These are the values of $\beta$ and $\gamma$ from our standard parameters---see Table \ref{table:baselineparams}.
}
\label{table:baselineparamsallcases}
\end{table}

% %%%%%%%%%%%%%%%%%%%%%%%%%%%%%%%%%%%%%%%%%%%%%%%%%%%%%%%%%%%%%%%%%%%%%%%%%%%%%
% Sens TABLE: OAT parameter sensitivity percentages
% %%%%%%%%%%%%%%%%%%%%%%%%%%%%%%%%%%%%%%%%%%%%%%%%%%%%%%%%%%%%%%%%%%%%%%%%%%%%%

\begin{table}[htbp]
\centering
  \begin{tabular}{c l p{3cm} p{3cm}}
    \hline\noalign{\smallskip} 
    Symbol & Description & $\Delta$ delay length & $\Delta$ dormancy length \\
    \noalign{\smallskip}\thickhline\noalign{\smallskip}
    
    $r_C$ & logistic growth rate of cancer & -50\% & 1\% \smallskip\\
    $r_{\text{max}}$ & maximum effective growth rate of cancer & 0\% & -7\%\smallskip\\
    $C^*$ & cancer steady state concentration & -83\% & -1\% \smallskip\\
    $\kappa$ & killing rate of cancer cells by T cells & +51\% & 0\% \smallskip\\
    
    \hline\noalign{\smallskip} 
    % ANTIGEN
    $r_A$ & antigen presentation source rate & -78\% & +1\%  \smallskip\\
    $\delta_A$ & antigen presentation degradation rate & 78\% & -1\%  \smallskip\\

    % INFLAMMATORY SIGNALS
    $r_I$ & inflammation source rate & -80\%  & +1\% \smallskip\\
    $\delta_I$ & inflammation degradation rate & 80\%  & -1\% \smallskip\\
    
    \hline\noalign{\smallskip}
    % T cells
    $r_E$ & effector cell growth coefficient & -74\%  & -1\% \smallskip\\
    $E^*$ & effector T cell base steady state & -74\% & -1\% \smallskip\\
    $r_S$ & non-effector T cell growth coefficient & -74\%  & -1\% \smallskip\\
    $S^*$ & non-effector T cell base steady state & 74\%  & 1\% \smallskip\\
    $\beta$ & effector T cell recruitment coefficient & -74\% & -4\% \smallskip\\
    $\gamma$ & non-effector T cell recruitment coefficient & 83\%  & +1\% \smallskip\\
    \hline\noalign{\smallskip}
  \end{tabular}\smallskip
  \caption{\bb{Parameter sensitivities.} Starting with the baseline parameters that simulate a 2-month delay (see Table \ref{table:baselineparams}), we deduce the change in the delay length and dormancy length by varying the value of the parameters one at a time by 1\%. The $+$ indicates a longer length of time, while $-$ indicates shorter. Overall, we can see that every parameter (with the exception of $r_\text{max}$) substantially affects the delay length, but not the dormancy length. On the other hand, $r_\text{max}$ is the only parameter that affects the dormancy delay to a nontrivial level.}
  \label{table:parameterOATsensitivity}
\end{table}

% %%%%%%%%%%%%%%%%%%%%%%%%%%%%%%%%%%%%%%%%%%%%%%%%%%%%%%%%%%%%%%%%%%%%%%%%%%%%
% Sens DIAGRAM: delay length vs parameter
% %%%%%%%%%%%%%%%%%%%%%%%%%%%%%%%%%%%%%%%%%%%%%%%%%%%%%%%%%%%%%%%%%%%%%%%%%%%%

%%%%%%%%%%%%%%%%%%%%%%%%%%%%%%%%%%%%%%%%%%%%%%%%%%%%%%%%%%%%%%%%%
%%%%%%%%%%%%%%%%%%%%%%%%%%%%%%%%%%%%%%%%%%%%%%%%%%%%%%%%%%%%%%%%%
%%%%%%%%%%%%%%%%%%%%%%%%%%%%%%%%%%%%%%%%%%%%%%%%%%%%%%%%%%%%%%%%%
%%%%%%%%%%%%%%%%%%%%%%%%%%%%%%%%%%%%%%%%%%%%%%%%%%%%%%%%%%%%%%%%%
%%%%%%%%%%%%%%%%%%%%%%%%%%%%%%%%%%%%%%%%%%%%%%%%%%%%%%%%%%%%%%%%%
%%%%%%%%% 6. DISCUSSION %%%%%%%%%%%%%%%%%%%%%%%%%%%%%%%%%%%%%%%%%
%%%%%%%%%%%%%%%%%%%%%%%%%%%%%%%%%%%%%%%%%%%%%%%%%%%%%%%%%%%%%%%%%
%%%%%%%%%%%%%%%%%%%%%%%%%%%%%%%%%%%%%%%%%%%%%%%%%%%%%%%%%%%%%%%%%
%%%%%%%%%%%%%%%%%%%%%%%%%%%%%%%%%%%%%%%%%%%%%%%%%%%%%%%%%%%%%%%%%
%%%%%%%%%%%%%%%%%%%%%%%%%%%%%%%%%%%%%%%%%%%%%%%%%%%%%%%%%%%%%%%%%
%%%%%%%%%%%%%%%%%%%%%%%%%%%%%%%%%%%%%%%%%%%%%%%%%%%%%%%%%%%%%%%%%

\newpage
\section{Discussion}\label{sec:discussions}

Our goal with this paper is to stimulate interest among math biologists in a largely ignored niche in the community, namely the phenomena of delayed responses seen in some patients. To do this, we introduced a foundational ODE model that balances a number of interests. Firstly, our model generates a variety of different response types observed experimentally, notably the interesting delayed case, without leveraging any artificial mechanisms to generate delays, such as delay differential equations. To highlight the primacy of our interest in the delayed-response case, the majority of our simulations assume a set of baseline parameters (see Table \ref{table:baselineparams}) that generate a two-month delay. Moreover, our simulations are designed to be idealised in that emphasis is placed on the essence of the phenomenon over the realistic fitting of datasets, which are largely limited. For a given patient, the administration of discrete dosages of checkpoint blockades have the effect of modifying these parameters by a discrete amount, $\Delta\beta$ or $\Delta\gamma$. An insight of our model is that the region within the ($\beta, \gamma$) parameter space generating delayed responses is very thin. This implies it could be relatively rare for a patient to land in the delayed response region after a dosage of the drug. Beyond immune profiles and drug efficacy, our model also offers insights on the balance between effector and non-effector T cells during the anti-tumour response. In particular, our model predicts that the effector and non-effector T cells ($E$ and $S$) cells share the tumour microenvironment on briefly, and one population dominates the other in between crucial transition points. We hope that our model, with its emphasis on simplicity and the high-level phenomenon, serves as a springboard for future work.

%%%%%%%%%% 6.1 Types of respoonses generated by model
\subsection{Four types of responses}

Four parameters in our model generate the spectrum of response types. They are $\beta$, $\gamma$, $E^*$, and $r_{\text{max}}$. The first and second parameters, $\beta$ and $\gamma$, control the respective levels of CTLA-4 and PD-1 expression and their therapeutic inhibitor in a given individual. Values for this pair of parameters can be calibrated in a logical manner to reflect individual inhibitors (e.g. anti-CTLA-4 Ipilimumab or anti-PD-1 Nivolumab) or their combination. A quicker anti-tumour response corresponds to increasing $\beta$ or decreasing $\gamma$. Focusing on delayed responses from our model, we note that the precise length of delay can be continuously calibrated up to a duration spanning years. Moreover, this can be obtained by perturbing $\beta$ and $\gamma$ alone. Our model also predicts that combining individual inhibitors provides therapeutic benefits. For instance, a CTLA-4 inhibitor with a 5-month delay combined with a PD-1 inhibitor with a 4-month delay results in a treatment with just 2-months delay. The synergistic effect of combination anti-CTLA-4 and PD-1 has been widely observed in trials~\cite{buchbinderCTLA4PD1Pathways2016a, rotteCombinationCTLA4PD12019}. Because $\beta$ and $\gamma$ alone are able to calibrate the delay length with great effect, our results suggest that amongst all possible factors that explain the existence of a delay for immune checkpoint inhibitor therapy, the average level of CTLA-4 and PD-1 expression on an individual alone may play key roles in determining their own response to this type of therapy. Moreover, the high-level scope at which we have defined $\beta$ and $\gamma$ in our model makes it amenable to existing and competing theories on CTLA-4's mechanism of action. One paradigm argues that anti-CTLA-4 removes barriers towards naive T cell activation, while another contends that the inhibitor increases anti-cancer activity by depleting regulatory T cells at the tumour site~\cite{rowshanravanCTLA4MovingTarget2018, tangAntiCTLA4AntibodiesCancer2018}. The former theory would imply we associate only $\beta$ to CTLA-4 expression, while the latter implies we might identify both $\beta$ and $\gamma$ as controllers of CTLA-4. In either case, our model is able to replicate delayed responses.\\
\\
\indent The third parameter, chosen to be $E^*$ in our baseline parameters, relates to patient-specific characteristics of an individual's immune system. However, we could also choose other parameters such as $S^*$, $r_A$ and $r_I$ in its place and maintain the ability for the model to generate the different response types. The parameter $E^*$ describes the natural strength of an individual's immune system against cancer, where higher values are associated with quicker anti-tumour responses and shorter delay durations. A particular interpretation could be $E^*$ measures the level of an individual's tumour-infiltrating lymphocytes. In contrast, $S^*$ describes the cancer's natural ability to suppress the host immune system, where higher values result in longer delay durations. The parameters $r_A$ and $r_I$ influence the rate at which tumour-specific T cells are mobilised by antigen and inflammation. Higher values of these parameters correlate with quicker anti-tumour responses and shorter delay lengths. All of these parameters describe patient-specific factors in some matter---whether it be at the antigen generation stage or the tumour-immune interaction stage.\\
\\
\indent Finally, the fourth parameter $r_{\text{max}}$ calibrates the maximum cancer growth rate. Setting a high ceiling for this rate allows us to obtain the quick and partial response, whereby the cancer will respond quickly to treatment but does not become eradicated. A limitation of the model is that it may be necessary to elevate $r_{\text{max}}$ to a biologically unrealistic level. Hence, the a quick partial response may be a rare occurrence in reality.

%%%%%%%%%% 6.2 Delays vs parameters
\subsection{Parameter insights}\label{subsec:delayvsparameter}

A set of observations can be made about model parameters and outputs, which subsequently lead to potential medical insights that we shall describe in Section \ref{subsec:immuneprofileanddrugefficacy}. Firstly, the delay length is sensitive to almost every parameter. Thus our model suggests that various patient-specific factors related to how an individual's immune system interacts with a tumour all play a part in determining the delay length from immune checkpoint inhibitors. These factors include the antigenicity of the tumour to the patient (controlled by $r_A$), the level of inflammatory response by the patient ($r_I$); and how the tumour affects and is affected by the patient's CTLA-4 and PD-1 expression levels ($\beta$ and $\gamma$). Patient-specific factors have been observed to affect treatment outcomes. For instance, candidate suitability for anti-PD-1 inhibitors are determined on the genetic profile of the tumour~\cite{hwangImmuneGeneSignatures2020, ayersIFNgRelatedMRNA2017}---factors that are specific to the individual.\\%

\indent Secondly, as we vary parameters, the delay length changes non-linearly. Moreover, it may only increase up to a point, upon which it becomes infinite. A diagrammatic explanation can be found in Fig.~\ref{SpikeODE_result_delay_vs_parameters_v4}. In short, the delay is sensitive to multiple parameters that control aspects of how a patient's immune system interacts with a tumour, and this relationship is non-linear and discontinuous. Importantly, this relationship points to the existence of a thin region within the parameter space where the delay exists, and slight perturbations can result in either the delay vanishing (quick response) or becoming infinite (no response).\\

\indent Thirdly, the model produces some counter-intuitive correlations between delay length and some parameters. Notably, the cancer-killing parameters and their effect on the delay length run counter to expectations. The model predicts that larger and faster growing tumours result in shorter delays, while patients with cancer cells that were killed quicker by effector T cells were nevertheless subject to longer delay lengths.

%%%%%%%%%% 6.3 Immune profiles and drug efficacy
\subsection{Immune profiles and drug efficacy}\label{subsec:immuneprofileanddrugefficacy}

Immunological differences between individuals can play a role in determining how they cope with cancer and react to treatments. For instance, key parameters like $E^*$ and $\gamma$, which relate to the level of tumour-infiltrating lymphocyte (TIL) and PD-1 expression in an individual, are relevant variables for treatment. For instance, higher levels of TIL are linked to better treatment outcomes~\cite{dushyanthenRelevanceTumorinfiltratingLymphocytes2015, paijensTumorinfiltratingLymphocytesImmunotherapy2021}, while suitability of checkpoint blockade treatment targeting PD-1 can depend on the genetic profile of the patient's tumour~\cite{hunterPDL1TestingGuiding2018}--these are all factors specific to individuals.\\

As a result, an individual immune profile can be considered as a set of parameter values. In particular, as we are mostly interested in the effect of checkpoint blockades, the primary parameters of interest are $\beta$ and $\gamma$. They broadly describe the activation and suppression rate of effector T cells, and serve as proxy variables for the levels of CTLA-4 and PD-1 expression. A CTLA-4 blockade has the effect of reducing CTLA-4 expression, which corresponds to an increase in the value of $\beta$ (higher Teff activation rate) for an individual. On a similar note, a PD-1 blockade has the effect of reducing PD-1 expression, which here, corresponds to a decrease in the value of $\gamma$ (lower Teff cell suppression rate) for an individual.\\

Zooming out, a population of patients map to different points in the parameter space of our model. These different immune profiles are distributed as different values of $\beta$ and $\gamma$ in Fig.~\ref{SpikeODE_result_delay_vs_parameters_v4}. In particular, patients yet to receive checkpoint blockades begin their treatment journey in the no-response parameter region. Patients who have not yet received any treatment lie in the no-response parameter region. The effect of administering immune checkpoint blockade therapy can then be modelled simply and phenomenologically by perturbing the value(s) of $\beta$ (for anti-CTLA-4), $\gamma$ (for anti-PD-1) or $(\beta, \gamma)$ (for a combination treatment). This allows us to describe the administration of discrete dosages, where each dosage would increase the patient's $\beta$ and decrease $\gamma$ with the step sizes $\Delta\beta$ and $\Delta\gamma$ commensurate to the dosage strength. Checkpoint blockade therapy is typically administered over months, with a minimum of a couple of weeks between dosages~\cite{rennerImmuneCheckpointInhibitor2019}.\\

\begin{figure}
    \includegraphics[scale=1]{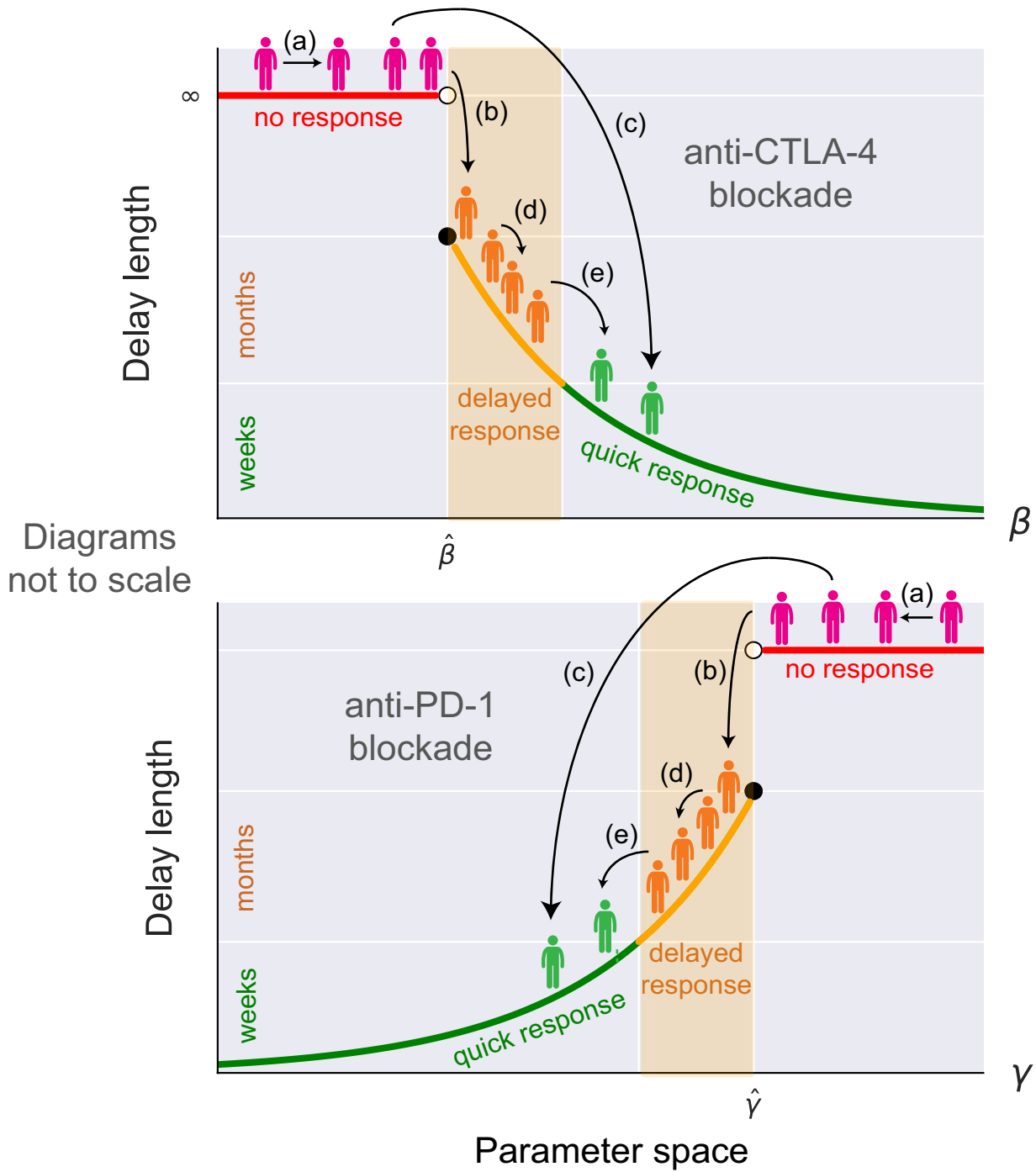}
    \renewcommand{\thefigure}{11}
    \caption{\bb{Illustrative diagram of relationship between immune profiles and checkpoint blockade treatment.} The key growth and suppression parameters $\beta$ and $\gamma$ broadly control the activation and suppression rate of T cells and vary between individuals. Thus the distribution of immune profiles across the population can be broadly modelled by a distribution spanning $\beta$ and $\gamma$. Moreover, $\beta$ and $\gamma$ serve as proxies for CTLA-4 and PD-1 expression in our model. For a given individual (which maps to a certain $\beta-\gamma$ combination), anti-CTLA-4 and anti-PD-1 blockades have the respective effects of increasing their personal value of $\beta$ and decreasing $\gamma$. This gives rise to a variety of cases. Blockades administered to an initial \ii{no-response patient} can be ineffective (a), have therapeutic effect albeit with a delay (b) or have high efficacy resulting in a quick response to treatment (c). Moreover, blockades administered to a \ii{delayed-response patient} can be ineffective (d) or sufficiently effective to overcome the delay (e).}
    \label{SpikeODE_result_delay_vs_parameters_drugs}
    \end{figure}

With this understanding of how our model describes treatment, the general idea is that a tumour-bearing individual start off with an initial ($\beta$, $\gamma$) that places them in the no-response parameter region, and successive dosages of CTLA-4 and/or PD-1 blockade move them away from this region in discrete steps towards the quick response region. This is diagrammatically summarised in Fig. \ref{SpikeODE_result_delay_vs_parameters_drugs}. At each dosage point, there may exist several possibilities. For a patient that falls in the no-response case after the $n$-th dosage (that is, $\beta < \hat{\beta}$ for anti-CTLA-4 and $\gamma > \hat{\gamma}$ for anti-PD-1), the $(n+1)$-th dosage gives three possibilities: still no-response, a response albeit delayed, or high drug efficacy resulting in a quick response. On the other hand, a patient that falls in the delayed response case after the $n$-th dosage ($\beta > \hat{\beta} + \epsilon$ for anti-CTLA-4 and $\gamma < \hat{\gamma} - \epsilon$ for anti-PD-1), the $(n+1)$-th dosage gives two possibilities: the response is still delayed, or the dosage was potent enough to upgrade them to the quick-response region.\\

The thinness of the parameter region that generates delayed responses may have interesting consequences. For instance, if we assume equal-strength dosages always manifest as equal step movements $\Delta\beta$ and/or $\Delta\gamma$ in the parameter space, then it is unlikely for a patient to land in the delayed-response region after a dose. In other words, our model predicts that the delayed responses are relatively rare. As a caveat, this assumes that immunological profiles across the population are distributed evenly across the parameter space, which might not be true. For instance, it might be that a significant proportion of people carry immune profiles that naturally place them in the delayed-response parameter region, even if it is mathematically thin.\\

These insights present opportunities for future work. Analyses can be done to try and see how the immune profiles of the population map to the parameter space and whether there are any clusters. Moreover, to help capture the inevitable variance between individual patients, our model could be extended by purposefully introducing noise in the parameters, so that variances between individual patients are explicitly modelled. For simplicity, the focus of this study is primarily on demonstrating a proof-of-concept backed by simulations.

% %%%%%%%%%%%%%%%%%%%%%%%%%%%%%%%%%%%%%%%%%%%%%%%%%%%%%%%%%%%%%%%%%%%%%%%%%%%%
% FIG 11. NEW Jun 2021. Immune profiles vs drugs
% %%%%%%%%%%%%%%%%%%%%%%%%%%%%%%%%%%%%%%%%%%%%%%%%%%%%%%%%%%%%%%%%%%%%%%%%%%%%

%%%%%%%%%% 6.4 T cell dynamics at the tumour site
\subsection{T cell insights}

The simulations offer several insights on T cell dynamics in the tumour microenvironment. Firstly, the dynamics between effector T cells, $E$, and non-effector T cells, $S$, can, from a modelling perspective, be described as a consensus-building model, like a competitive race between effector T cells (T eff) and regulatory T cells (T reg). In our simulations, the cancer population remains high due to a high level of $S$ cells and a simultaneously low level of $E$ cells. In the delayed response outcome, the eventual anti-tumour response is carried by a sudden influx of $E$ cells that is accompanied by the sudden exodus of $S$ cells. It is as if, after an initial transient period, the balance between $E$ and $S$ switches and now the $E$ population temporarily dominates the tumour site. The anti-tumour $E$ response is never durable, and their transient rise is accompanied by a return of both $S$ cells and the cancer, which starts the cycle anew. Thus the T cell populations behaves in an ongoing dynamical dance with momentum swinging quickly and decisively between $S$ and $E$ or vice versa. Moreover, this competition between $E$ and $S$ cells at the tumour site draws a parallel with the competition between effector T cells and regulatory T cells respectively, where the presence of effector T cells reduces the cancer population while regulatory T cells serve to maintain immune tolerance of the tumour and is associated with a poor prognosis~\cite{oleinikaSuppressionSubversionEscape2013, quezadaCTLA4BlockadeGMCSF2006, curranTumorVaccinesExpressing2009}.\\
\\
\indent Secondly, the model predicts that although T cell responses (where $E \approx E^*$ and $S \approx 0$) can be relatively short-lived, that response can hold off the cancer from returning for a much longer period. For instance, Fig.~\ref{SpikeODE_result_T_cell_responses_long} displays a prototypical situation for our model where an $E$ cell response that lasts for mere days is able to suppress the return of cancer for months.\\
\\
\indent A limitation of our model is that we do not model the level of CTLA-4 and PD-1 expression on individual T cells. For simplicity, we fix their average level of expression at a population level (using the parameters $\beta$ and $\gamma$), whereby a higher level represents no treatment and a lower level represents the administration of antibody inhibitors. For future work, it would be instructive to model individual T cells having varying levels of CTLA-4 and PD-1 expressions. Such enhancements would allow us to explicitly model naive and effector T cells, thereby providing the ability to test existing theories, such as how CTLA-4 affects the activation pathway of naive T cells and how PD-1 affects an individual T cell's cytotoxicity towards cancer cells. These extensions could be implemented by using a more extensive set of ODEs, a PDE model, or agent-based modelling techniques.

%%%%%%%%%%%%%%%%%%%%%%%%%%%%%%%%%%%%%%%%%%%%%%%%%%%%%%%%%%%%%%%%%%
%%%%%%%%%%%%%%%%%%%%%%%%%%%%%%%%%%%%%%%%%%%%%%%%%%%%%%%%%%%%%%%%%%
%%%%%%%%%%%%%%%%%%%%%%%%%%%%%%%%%%%%%%%%%%%%%%%%%%%%%%%%%%%%%%%%%%
%%%%%%%%%%%%%%%%%%%%%%%%%%%%%%%%%%%%%%%%%%%%%%%%%%%%%%%%%%%%%%%%%%
%%%%%%%%%%%%%%%%%%%%%%%%%%%%%%%%%%%%%%%%%%%%%%%%%%%%%%%%%%%%%%%%%%
%%%%%%%%% 7. BIBLIOGRAPHY %%%%%%%%%%%%%%%%%%%%%%%%%%%%%%%%%%%%%%%%
%%%%%%%%%%%%%%%%%%%%%%%%%%%%%%%%%%%%%%%%%%%%%%%%%%%%%%%%%%%%%%%%%%
%%%%%%%%%%%%%%%%%%%%%%%%%%%%%%%%%%%%%%%%%%%%%%%%%%%%%%%%%%%%%%%%%%
%%%%%%%%%%%%%%%%%%%%%%%%%%%%%%%%%%%%%%%%%%%%%%%%%%%%%%%%%%%%%%%%%%
%%%%%%%%%%%%%%%%%%%%%%%%%%%%%%%%%%%%%%%%%%%%%%%%%%%%%%%%%%%%%%%%%%
%%%%%%%%%%%%%%%%%%%%%%%%%%%%%%%%%%%%%%%%%%%%%%%%%%%%%%%%%%%%%%%%%%

\newpage

\bibliographystyle{plain} %spbasic}      % basic style, author-year citations
\bibliography{references.bib}   % name your BibTeX data base



% \begin{thebibliography}{}
% %
% % and use \bibitem to create references. Consult the Instructions
% % for authors for reference list style.
% % Format for Journal Reference
% Author, Article title, Journal, Volume, page numbers (year)
% Author, Book title, page numbers. Publisher, place (year)
% \end{thebibliography}

\end{document}